\documentclass[a4paper]{article}

\usepackage[english]{babel}

\usepackage[utf8]{inputenc}
\usepackage[utf8]{inputenc}
\usepackage{amsmath}
\usepackage{graphicx}
\usepackage{amssymb}
\usepackage{amsthm}
\usepackage{tikz-cd}
\usepackage{mathrsfs}
\usepackage[colorinlistoftodos]{todonotes}
\usepackage{enumitem}
\usepackage{yfonts}
\usepackage{ dsfont }
\usepackage{MnSymbol}
\usepackage{slashed}

\title{SYZ conjecture for Calabi-Yau hypersurfaces in the Fermat family}

\author{Yang Li}

\date{\today}
\newtheorem{thm}{Theorem}[section]
\newtheorem{lem}[thm]{Lemma}

\theoremstyle{definition}
\newtheorem{eg}[thm]{Example}

\newtheorem{cor}[thm]{Corollary}

\newtheorem{rmk}[thm]{Remark}
\newtheorem{prop}[thm]{Proposition}
\newtheorem{Def}[thm]{Definition}

\newtheorem*{Notation}{Notation}

\newtheorem*{Acknowledgement}{Acknowledgement}

\newcommand{\ie}{\emph{i.e.} }
\newcommand{\cf}{\emph{cf.} }
\newcommand{\into}{\hookrightarrow}

\newcommand{\R}{\mathbb{R}}
\newcommand{\C}{\mathbb{C}}
\newcommand{\Z}{\mathbb{Z}}

\newcommand{\Q}{\mathbb{Q}}

\newcommand{\norm}[1]{\left\lVert#1\right\rVert}
\newcommand{\Lap}{\Delta}

\DeclareMathOperator{\Hom}{Hom}

\def\Xint#1{\mathchoice
	{\XXint\displaystyle\textstyle{#1}}%
	{\XXint\textstyle\scriptstyle{#1}}%
	{\XXint\scriptstyle\scriptscriptstyle{#1}}%
	{\XXint\scriptscriptstyle\scriptscriptstyle{#1}}%
	\!\int}
\def\XXint#1#2#3{{\setbox0=\hbox{$#1{#2#3}{\int}$ }
		\vcenter{\hbox{$#2#3$ }}\kern-.6\wd0}}

\def\dashint{\Xint-}

\begin{document}
	\maketitle

\begin{abstract}
We produce special Lagrangian $T^n$-fibrations on the generic regions of some Calabi-Yau hypersurfaces
in the Fermat family
$X_s=\{ Z_0\ldots Z_{n+1}+ e^{-s} (  Z_0^{n+2}+ \ldots Z_{n+1}^{n+2} ) =0   \}\subset \mathbb{CP}^{n+1}
$ near the large complex structure limit $s\to +\infty$.
\end{abstract}

\section{Introduction}

The Strominger-Yau-Zaslow (SYZ) conjecture \cite{SYZ} is the following: given a family of $n$-dimensional polarised Calabi-Yau (CY) manifolds $(X_s, g_s, J_s, \omega_s, \Omega_s)$ of holonomy $SU(n)$ degenerating to the large complex structure limit, then
\begin{itemize}
\item After suitable scaling, the metric spaces $(X_s, g_s)$ converge in the Gromov-Hausdorff sense to a \emph{singular affine manifold} $B$ homeomorphic to $S^n$. The limiting metric $g_\infty$ is a \emph{real Monge-Amp\`ere metric} on the smooth locus. (This part is also known as the Kontsevich-Soibelman conjecture \cite{KS}\cite{KS2}.)
\item
Near the degenerating limit, the manifold $X_s$ admits a \emph{special  Lagrangian $T^n$ fibration}  over the base $B$ with some singular fibres. The diameters of the fibres are much smaller compared to $\text{diam}(B)$. In the generic region on $X_s$, which covers most of the measure on $X_s$, the metric $g_s$ is a small perturbation of a \emph{semiflat metric}, meaning that the $T^n$ fibres are almost flat.
\item
Mirror manifolds should be constructed as another $T^n$ fibration over the same base $B$, by fibrewise replacing the $T^n$ fibres with the dual tori.
\end{itemize}

An early achievement is Gross and Wilson's gluing construction \cite{GrossWilson} of degenerating CY metrics on K3 surfaces with elliptic fibrations, which becomes a special Lagrangian $T^2$-fibration after hyperk\"ahler rotation. In this setting the metric is known semi-explicitly. The same period brought forth many insights concerning topological \cite{Gross2}, combinatorial \cite{HaaseZharkov}\cite{HaaseZharkov2}, and differential geometric \cite{Zharkov} aspects of the SYZ conjecture, until Joyce \cite{Joyce} discovered through his study of special Lagrangian singularities that the SYZ fibration map cannot be na\"ively expected to be smooth, indicating the difficulty of the metric problem.

Later research on the SYZ conjecture gradually shifted focus  from its metric geometric roots, in favour of softer approaches based on algebraic or symplectic methods, taking the original SYZ conjecture mainly as an inspiration. This has led to spectacular progress in the mathematical understanding of mirror symmetry, described in the excellent survey \cite{Gross}.

In the metric vein, the SYZ conjecture fits into the more general question of understanding how CY metrics degenerate as the complex and K\"ahler structures vary. The main dichotomy is whether the family of metrics are \emph{noncollapsed}, meaning there is a uniform lower bound on the volume once the diameter is normalised to one. In the noncollapsing case much is known: for example, a polarised family of noncollapsed CY manifolds can only degenerate to normal CY varieties with klt singularities, and the notion of metric convergence agrees with the algebro-geometric notion of flat limit  \cite{DonaldsonSun}.

The collapsing case is widely open. Tosatti et al. made substantial progress on describing collapsing metrics associated with holomorphic fibrations \cite{Tosatti}\cite{GrossTosattiZhang}, in particular generalising much of \cite{GrossWilson} to hyperk\"ahler manifolds with holomorphic Lagrangian fibrations. Recently there are many efforts to describe the degenerating CY metrics in special cases, notably for K3 surfaces \cite{Lorenzo}\cite{HSVZ}\cite{Odaka}, and higher dimensional generalisations \cite{HSVZ2}.

The metric SYZ conjecture resisted most attempts because the large complex structure limit is a very severe degeneration mechanism. An interesting program of Boucksom et al. \cite{Boucksom}\cite{Boucksom1} proposes that in the case of polarised algebraic degenerations the underlying Calabi-Yau manifolds  converge naturally into a non-archimedean (NA) space, and the CY metrics should converge in a potential theoretic sense to their NA analogue. Their greatest achievements so far is to define and solve the NA Monge-Amp\`ere (MA) equation, building on heavy machinery from birational geometry. To make contact with the SYZ conjecture, it would still remain to compare the non-archimedean MA equation with the real MA equation, prove the potential theoretic convergence, and improve it to the metric convergence. Notwithstanding these difficulties, this program has the promise to prove the SYZ conjecture in great generality. 

The viewpoint of this paper is much more concrete.
We focus on the Fermat family of projective hypersurfaces of any dimension $n$, approaching the large complex structure limit:
\begin{equation}
X_s= \{  Z_0Z_1\ldots Z_{n+1}+ e^{-s} \sum_{i=0}^{n+1} Z_i^{n+2}=0     \}, \quad s\gg 1.
\end{equation}
The most striking aspect of our results is
\begin{thm}\label{mainthm}
	(\cf section \ref{SpecialLagrangianfibraiongenericregion})
For the Fermat family, consider the Calabi-Yau metrics on $X_s$ in the polarisation class $s^{-1} [\Delta]$ where $[\Delta]$ is a fixed K\"ahler class on $\mathbb{CP}^{n+1}$ restricted to $X_s$. Then for a subsequence of $X_s$ as $s\to +\infty$, there exists a special Lagrangian $T^n$-fibration on the generic region $U_s\subset X_s$, such that $\frac{\text{Vol}(U_s) }{ \text{Vol}(X_s) }\to 1$ as $s\to +\infty$.
\end{thm}

We also summarize informally the other results in this paper:
\begin{itemize}
\item  (\cf section \ref{GromovHausdorffconvergencesection})  The subsequence of CY metrics converge in the Gromov-Hausdorff sense to the metric completion of a smooth real MA metric on an open dense subset $\mathcal{R}\subset \partial \Delta_\lambda^\vee$, where $\partial \Delta_\lambda^\vee$ denotes the boundary of a certain $(n+1)$-dimensional simplex $ \Delta_\lambda^\vee$ in $\R^{n+1}$ arising naturally from tropical geometry, and $\partial \Delta_\lambda^\vee\setminus \mathcal{R}$ has zero $(n-1)$-Hausdorff measure.

\item   (\cf Prop. \ref{diameterbound})  The diameters of the subsequence of CY metrics are uniformly bounded.

\item (\cf section \ref{Higherregularitysection}) In the generic region of $X_s$ for $s\gg 1$, the CY metrics are $C^\infty_{loc}$ close to a sequence of semiflat metrics. In particular the sectional curvature in the generic region is uniformly bounded.
\end{itemize}

A basic feature of the complex geometry of CY hypersurfaces near the large complex structure limit, is that in generic regions  the local structure is a large annulus region in $(\C^*)^n$, equipped with a holomorphic volume form which modulo a scale factor is very close to $d\log z_1\wedge \ldots d\log z_n$.
An elementary observation is that \emph{plurisubharmonic (psh) functions are intimately related to convex functions}:
\begin{itemize}
\item  Let $\phi$ be psh on an annulus $\{ 1< |z_j|< \Lambda  \}\subset (\C^*)^n$, then the fibrewise average function
\[
\bar{\phi}(x_1,\ldots x_n)= \dashint_{T^n} \phi( e^{x_1+ i\theta_1}, \ldots e^{x_n+i\theta_n}  ) d\theta_1\ldots d\theta_n
\]
is convex.

\item
Let $u$ be a convex function on $\{ 0<x_j< \log \Lambda \}$, then the pullback of $u$ to $\{ 1< |z_j|< \Lambda  \}\subset (\C^*)^n$ via the logarithm map is psh, and $u$ solves the real MA equation $\det(D^2 u)=\text{const}$ iff its pullback solves the complex MA equation $\det ( \frac{\partial^2 u}{\partial \log z_i \partial \overline{ \log z_j } }   )= \text{const}$.

\end{itemize}

Our strategy is to show that in the highly collapsed regime $s\gg1$, the local K\"ahler potentials are $C^0$-approximated by convex functions, whose regularity properties can be then transferred back to the local K\"ahler potentials at least in the generic region. In effect, this implies in the generic region the Calabi-Yau metrics are collapsing with uniformly bounded sectional curvature; then the existence of the special Lagrangian fibration in the generic region is a simple perturbation argument. Keeping in mind that the local complex structure is an annulus in $(\C^*)^n$, the special Lagrangian fibration is just a small $C^\infty$-perturbation of the logarithm map $(\C^*)^n\to \R^n$.

 The essential problem is to obtain uniform estimates on the CY metrics as $s\to \infty$. Our techniques differ very significantly from Yau's proof of the Calabi conjecture. Our  K\"ahler potential estimates are largely based on Kolodziej's method in pluripotential theory, which has the advantage of robustness even in collapsing settings. The technical core of our contribution  is to produce a regularisation of the Calabi-Yau potential, and prove an improved version of the global Skoda inequality, which for large $s$ forces the potential to be very close to its regularisation. As convexity is built into the construction of the regularisation, this furnishes a bridge between holomorphic and convex geometry, and one can start to transfer the a priori much better regularity from the convex world into the holomorphic world near the collapsing limit $s\to \infty$.
 Our higher order estimates exploit the local regularity theory of real MA equations, and a result of Savin from nonlinear PDE theory.

The structure of the paper is as follows.
We survey the rather extensive analytical backgrounds in section  \ref{analyticalbackgrounds}. The complex geometry of the degenerating hypersurfaces is discussed in section \ref{Degeneratingtorichypersurfaces}, with particular emphasis on its interplay with tropical geometry. We estimate the Calabi-Yau potentials in section \ref{Estimatesonpotential}; in particular we prove the Skoda type estimates, the uniform $L^\infty$ bound, and the $C^0$-approximation by the convex regularisations. In section \ref{FermatcasemetricconvergenceandSYZfibration}, we use uniform Lipschitz bounds on the regularisation to extract a subsequential limit, and show that this defines a real MA metric. We then use the local regularity theory of real MA metrics to show the higher order estimates on the CY local potentials, and prove the existence of the special Lagrangian fibration.

We now discuss some directions of future research. 
\begin{itemize}
\item
It seems highly plausible that the SYZ conjecture on generic regions will hold also on many other degenerating CY manifolds, or at least CY hypersurfaces. In fact the only reason we restrict to the Fermat case is to utilize the large discrete symmetry group to give a relatively simple proof of a technical extension property for locally convex functions, which seems likely to generalise to other contexts. 

\item 
One would like to study the existence, uniqueness, and regularity of the real MA equation on compact polyhedral sets, which are covered by charts whose transition functions are only piecewise linear but not smooth in general; the SYZ conjecture predicts the solutions to such real MA equations should arise as possible limits of the collapsing CY metrics.  This question may be parallel to the non-archimedean MA approach taken up in \cite{Boucksomsurvey}. At present according to the author's knowledge, it is not clear how to define the real MA equation globally on such sets, and in fact we do not even have an established notion of local convexity.

Such questions on the real MA equations have direct bearings on improving our main theorem. For instance, if one can establish uniquenss, then there is no need to pass to subsequences in all of our results. If one can establish sufficient regularity, then it may be possible to prove the Gromov-Hausdorff limit is homeomorphic to $\partial \Delta_\lambda^\vee\simeq S^n$.

The problem to set up the real MA equation is quite subtle. On a piecewise linear manifold the notion of a convex function is dependent on charts, and so does the real MA operator. To set up an invariant notion of the real MA equation, it is necessary to make branch cuts to charts. The location of such cuts seems to depend on some gradient condition on the convex function in question, and is hard to predict in the absence of symmetry. Thus the global real MA equation on polyhedral sets has the feature of a free boundary problem.

\item
The a priori estimate approach in this paper says very little about the CY metrics in regions with high curvature concentration. In the case of CY 3-folds, the author  \cite{Li} recently constructed the 3-dimensional analogues of the Ooguri-Vafa metric, which are conjectured to be the universal metric models for the neighbourhood of the most singular fibres in a generic SYZ fibration. A program to tackle the 3-fold case of the SYZ conjecture based on gluing ideas is outlined in \cite{Li}, which has the ultimate aim to give a global description of the metric, and to produce a special Lagrangian fibration globally. This gluing approach requires very refined information on the singularities of the real MA equation, which is still far from what we can establish by a priori estimate considerations.

\end{itemize}

\begin{Acknowledgement}
The author is a postdoc at the IAS, funded by the Zurich
Insurance Company Membership. The pluripotential theoretic approach is inspired by the talks of Boucksom. The author would like to thank S. Sun, S.  Donaldson, Y. Jhaveri, C. Mooney and P. Sarnak for discussions, W. Feldman for giving a simple proof to a technical lemma, and the IAS for providing a stimulating research environment.
\end{Acknowledgement}

\section{Analytic backgrounds}\label{analyticalbackgrounds}


\subsection{Skoda inequality}

An upper semicontinuous $L^1$-function $\phi$ on a coordinate ball is called plurisubharmonic (psh)  if $\sqrt{-1}\partial \bar{\partial} \phi\geq 0$. The basic intuition is that regularity properties for psh functions in general dimensions  are analogous to subharmonic functions on Riemann surfaces. This is captured by the basic version of the \emph{Skoda inequality}:

\begin{thm}\label{Skodabasicversion}
(\cf \cite[Thm 3.1]{Zeriahi}) If $\phi$ is psh on $B_2\subset \C^n$, with $\int_{B_2} |\phi| \omega_E^n \leq 1$ with respect to the standard Euclidean metric $\omega_E$, then
there are dimensional constants $\alpha$, $C$, such that 
\[
\log \int_{B_1} e^{-\alpha \phi} \omega_E^n \leq  C.
\]

\end{thm}

\begin{rmk}
If instead $\int_{B_2} |\phi| \omega_E^n \leq C'$ for some constant $C'$, then we can apply Thm \ref{Skodabasicversion} to a scaling of $\phi$, to get a Skoda inequality with modified $\alpha, C$.
\end{rmk}

\begin{rmk}\label{ChernLevine}
	Assuming an $L^1$-bound on $\phi$, then 
	we can take a suitable cutoff function $\chi$, and via integration by parts,
	\[
	\int_{B_1} \sqrt{-1}\partial \bar{\partial} \phi\wedge \omega_E^{n-1} \leq   \int_{B_2}\chi \sqrt{-1}\partial \bar{\partial} \phi\wedge \omega_E^{n-1} = \int_{B_2} \phi \sqrt{-1}\partial \bar{\partial} \chi\wedge \omega_E^{n-1} \leq \norm{\chi}_{C^2} \norm{\phi}_{L^1} \leq C . 
	\]
	This simple idea is a basic version of the \emph{Chern-Levine inequality}, which is another fundamental reason why psh functions are much more regular than the subharmonic functions in general dimensions.
\end{rmk}

The basic Skoda inequality immediately implies a global version. On a compact K\"ahler manifold $(X, \omega)$, we say an upper semicontinuous $L^1$-function $\phi\in PSH(X,\omega)$  if $\omega_\phi=\omega+ \sqrt{-1} \partial \bar{\partial}\phi\geq 0$. This is the generalised notion of K\"ahler potentials.

\begin{thm}\label{Skodastandardversion}
On a fixed $(X, \omega)$, there are positive constants $\alpha$, $C$ depending only on $X, \omega$, such that
\[
\int_X e^{-\alpha \phi} \omega_X^n \leq C, \quad \forall \phi \in PSH(X,\omega) \text{ with } \sup \phi=0.
\]
\end{thm}

\begin{rmk}
Here $\int_X |\phi| \omega_X^n$ is automatically bounded  using the Harnak inequality, because $\Lap \phi \geq -n$ for $\phi\in PSH(X,\omega)$.
\end{rmk}

\begin{rmk}
The supremum of all such $\alpha$ is known as Tian's alpha invariant. 
\end{rmk}

\subsection{Kolodziej's estimate on pluripotentials}\label{Toolsfrompsh}

Here we outline a method to estimate K\"ahler potentials, pioneered by  Kolodziej, and further developed by \cite{DemaillyPali} and  \cite{EGZ}\cite{GZ}. Our exposition largely adapts  \cite{EGZ}\cite{GZ}\cite{EGZ2}, with special attention to the dependence of constants. Unlike in \cite{EGZ}, we do not impose a volume normalisation.

Given an $n$-dimensional K\"ahler manifold $(X,\omega)$, for $\phi\in PSH(X,\omega)\cap L^\infty$, pluripotential theory allows one to make sense of the Monge-Amp\`ere (MA) measure $\omega_\phi^n$, generalising the notion of volume forms. The basic problem is to estimate $\phi$  from a priori bounds on $\omega_\phi^n$.  A key concept is the \emph{capacity} of subsets $K\subset X$:
\[
Cap_\omega(K)= \sup \{  \int_K \omega_u^n | u\in PSH(X,\omega),  0\leq u\leq 1        \}.
\]

We wish to sketch the main ideas behind a prototypical result:

\begin{thm}\label{pluripotentialthm1}
	Let $(X, \omega)$ be a compact K\"ahler manifold, and $\phi\in PSH(X,\omega)\cap C^0$, such that $\omega_\phi^n$ is an absolutely continuous measure. Assume there are positive constants $\alpha, A$, such that the Skoda type estimate holds with respect to $\omega_\phi^n$:
	\begin{equation}\label{Skodaassumption}
	\int_X e^{-\alpha u}  \frac{\omega_\phi^n }{ \text{Vol}(X) } \leq A, \quad \forall u\in PSH(X,\omega) \text{ with } \sup_X u=0.
	\end{equation}
\begin{itemize}

\item   For fixed $n,\alpha,A$, there is number $B(n,\alpha,A)$, such that if $\frac{ \int_{ \phi\leq -t_0} \omega_\phi^n}{  Vol(X)  }< (2B)^{-2n}$ for some $t_0$, then $\min \phi\geq -t_0- 4B (\frac{ \int_{ \phi\leq -t_0} \omega_\phi^n}{  Vol(X)} )^{1/2n}$.  

\item   If $\sup_X \phi=0$, then $\norm{\phi}_{C^0} \leq C(n,\alpha,A)$.
\end{itemize}

\end{thm}


The first ingredient is:

\begin{lem}(\cf \cite[Lemma 2.3]{EGZ})\label{volumecontrolscapacity}
	The MA measure of sublevel sets controls the capacity of lower sublevel sets : for $\tau\geq 0$ and $0\leq t\leq 1$, 
	\[
	t^n Cap_\omega ( \phi<-\tau-t  ) \leq \int_{\phi<-\tau} \omega_\phi^n.
	\]
\end{lem}

The second ingredient below contains the most substance:

\begin{lem}
	(Volume-capacity estimate) In the setting of Thm. \ref{pluripotentialthm1}, for any compact set $K\subset X$,
	\begin{equation}\label{volumecapacity}
	\frac{\int_K \omega_\phi^n}{\text{Vol}(X)} \leq A e^\alpha \exp \left( -\alpha (\frac{ \text{Vol}(X)}{ \text{Cap}_\omega(K)})^{1/n}  \right) .
	\end{equation}
	In particular there is a constant $B=B(n,\alpha, A)$ verifying the power law bound
	\[
	\frac{\int_K \omega_\phi^n}{\text{Vol}(X)} \leq B^{2n} \left( \frac{ \text{Cap}_\omega(K)}{ \text{Vol}(X)}   \right)^2  .
	\]
\end{lem}

\begin{proof}
	(Sketch) 
	We may assume $K$ is not pluripolar, for otherwise $\int_K \omega_\phi^n=0$ and $\text{Cap}_\omega(K)=0$.
	We introduce the Siciak extremal function 
	\[
	V_{K,\omega}= \sup\{  u\in PSH(X, \omega)|u\leq 0 \text{ on } K      \},
	\]
	whose upper semicontinuous regularisation $V_{K,\omega}^*\in PSH(X,\omega)$. By the Alexander-Taylor comparison principle (\cf \cite[Prop. 6.1]{GZ}), 
	\[
	\exp(- \sup_X V_{K,\omega}) \leq e \exp\left( - (\frac{ \text{Vol}(X)}{ \text{Cap}_\omega(K)})^{1/n}  \right).
	\]
	By the Skoda integrability assumption (\ref{Skodaassumption}), and the fact that $V_{K,\omega}=V_{k,\omega}^*$ a.e with respect to $\omega^n$ (so by absolute continuity also for $\omega_\phi^n$),
	\[
	\int_X e^{  \alpha( \sup_X V_{K,\omega}- V_{K,\omega}    )} \omega_\phi^n =
	\int_X e^{  \alpha( \sup_X V_{K,\omega}- V_{K,\omega}^*    )} \omega_\phi^n \leq A\text{Vol}(X) ,
	\]
	hence
	\[
	\int_K e^{  -\alpha V_{K,\omega}   } \omega_\phi^n\leq
	\int_X e^{  -\alpha V_{K,\omega}    } \omega_\phi^n \leq A e^\alpha \text{Vol}(X) \exp \left( -\alpha (\frac{ \text{Vol}(X)}{ \text{Cap}_\omega(K)})^{1/n}  \right) .
	\]
	The volume-capacity estimate (\ref{volumecapacity}) follows because $V_{K,\omega}\leq 0$ on $K$.
\end{proof}

The third ingredient is an elementary decay lemma:

\begin{lem}
	(\cf \cite[Lemma 2.4 and Remark 2.5]{EGZ}) Let $f: [t_0,\infty)\to [0,\infty)$ be a nonincreasing right-continuous function, such that
	\[
	\begin{cases}
	f(t_0)< \frac{1}{2B}, 
	\\
	tf(\tau+t) \leq B f(\tau)^2 ,\quad \forall \tau\geq 0, \quad 0\leq t\leq 1,
	\\
	\lim_{t\to \infty} f(t)=0.
	\end{cases}
	\]
	Then $f(t)=0$ for $t\geq t_0+ 4B f(t_0)$.
\end{lem}

\begin{proof}
(Thm \ref{pluripotentialthm1})
 Combining the first two ingredients, the function $f(t)= (\frac{ \int_{\phi\leq-t} \omega_\phi^n }{\text{Vol}(X)  } )^{1/2n} $ satisfies
\[
t f(t+\tau) \leq  Bf(\tau)^2, \quad 0\leq t\leq 1, \quad \tau\geq 0,
\]
We conclude that for $t> t_0+ 4Bf(t_0)$ the sublevel set $\{ \phi\leq -t \}$ has zero $\omega_\phi$-measure, and therefore zero capacity by Lemma \ref{volumecontrolscapacity}, so $\phi$ has the lower estimate as claimed in the first statement.

For the second statement, by (\ref{Skodaassumption}) we have an a priori exponential decay
\[
f(t)\leq A^{1/2n} e^{-\alpha t/2n}, \quad t\geq 0,
\]
which allows us to find an appropriate $t_0$.	
\end{proof}

\begin{rmk}
	Thm. \ref{pluripotentialthm1} implies a famous result of Kolodziej stating that if we fix $(X,\omega)$ and $p>1$, then $\phi$ has a $C^0$-bound depending only on $X,\omega,\norm{ \frac{\omega_\phi^n}{\omega^n} }_{L^p}$. It is enough to check (\ref{Skodaassumption}), which reduces by H\"older inequality to the standard Skoda inequality (\cf Thm \ref{Skodastandardversion}), with modified constants.
	The strength of Thm. \ref{pluripotentialthm1} is that it still applies when the complex/K\"ahler structures are highly degenerate, as it distills the dependence on $(X,\omega)$ to only 3 constants $n,\alpha,A$.
\end{rmk}

Thm. \ref{pluripotentialthm1} gives a criterion for two K\"ahler potentials to be close to each other.

\begin{cor}\label{StabilityestimateKolodziej}
	(Stability estimate) Let $(X,\omega)$ be a compact K\"ahler manifold, and $\phi, \psi\in PSH(X,\omega)\cap C^0$, such that $\omega_\phi^n$ is absolutely continuous. Assume $\norm{\psi}_{C^0} \leq A'$ and the Skoda type estimate (\ref{Skodaassumption}). Then there is a number $B(n,A, A', \alpha)$, such that if $\frac{ \int_{ \phi-\psi\leq -t_0} \omega_\phi^n}{  Vol(X)  }< (2B)^{-2n}$ for some $t_0$, then \[
	\min (\phi-\psi)\geq -t_0- 4B \left(\frac{ \int_{ \phi-\psi\leq -t_0} \omega_\phi^n}{  Vol(X)} \right)^{1/2n}.
	\].   
\end{cor}

\begin{proof}
If $\psi$ is smooth, this follows from Thm. \ref{pluripotentialthm1} by changing $\omega$ into $\omega_\psi$, and changing $\phi$ into $\phi-\psi$, and checking the Skoda type estimate holds with modified constants. In general, one can approximate $\psi\in PSH(X,\omega)$ by a decreasing sequence of functions in $PSH(X,\omega)\cap C^\infty$ \cite{Blocki}, and since $\psi\in C^0$ the convergence is uniform by Dini's theorem. 
\end{proof}

\subsection{Algebraic metrics and asymptotes}\label{Algebraicmetrics}

This section is included for motivational purposes. On any compact complex manifold $X$ with a positive line bundle $L$, any fixed K\"ahler metric $\omega$ in the class $2\pi c_1(L)$ is the curvature form of a Hermitian metric $h$ on $L$. Consider the projective embedding  $\iota_k:X\into \mathbb{P}(H^0(X,L^k) ^*)$ for $k\gg 1$. The $L^2$ norms on sections induce Euclidean metrics on the vector spaces $H^0(X,L^k)$, hence Fubini-Study metrics $\omega_{FS,k}$ on $\mathbb{P}(H^0(X,L^k) ^*)$. A 
famous result of Tian says that $\omega$ is approximated by the algebraic metrics $k^{-1}\iota_k^*\omega_{FS,k}$ as $k\to \infty$; this idea has been much exploited in \emph{regularization theorems}.

This construction is particularly transparent in the toric case, as explained in \cite{Donaldsontoric}. Let $(X,L)$ be an $n$-dimensional polarised toric manifold with moment polytope $P$, so a $T^n$-invariant basis $\{ s_m \}$ of $H^0(X,L^k)$ corresponds to $kP\cap \Z^n$, or equivalently $P\cap k^{-1}\Z^m$ after rescaling. The $L^2$-metric on $H^0(X,L^k)$ is \emph{diagonal} in the basis; \ie the toric assumption reduces the unitary group acting on $H^0(X,L^k)$ to its maximal torus. Concretely, let $\phi$ denote the torus invariant K\"ahler potential on $X\cap (\C^*)^n$, equivalently thought as some convex function of $\vec{t}\in\R^n$ via the logarithm map $(\C^*)^n\to \R^n$. Then 
\begin{equation}\label{L2metrictoriccase}
I_m(k)=\norm{s_m}_{L^2}^2=\int_X |s_m|^2 d\text{Vol}= \text{const} \int_{\R^n} e^{-k(\phi- \vec{t}\cdot m  )} d\vec{t},\quad m\in P\cap k^{-1}\Z^n,
\end{equation}
and the Fubini-Study potentials are
\begin{equation}\label{FubiniStudypotential}
k^{-1}\iota_k^* \phi_{FS,k}= k^{-1} \log \left(   \sum_{m\in P\cap k^{-1}\Z^n} I_m(k)^{-1}  |s_m|^2 \right).
\end{equation}
Now the RHS of (\ref{L2metrictoriccase}) is a \emph{Laplace type integral}, and its dominant contribution comes from the neighbourhood of the point $\vec{t}_0$ where $\vec{t}\cdot m- \phi(\vec{t})$ is maximized among $\vec{t}\in \R^n$. The maximum is the value of the \emph{Legendre transform} of $\phi$:
\[
u(m)= \sup_{\vec{t}} (\vec{t}\cdot m- \phi(t) ).
\]
The steepest descent method yields the asymptote
\[
k^{-1} \log I_m(k)= u(m) +O( k^{-1}\log k ), \quad k\to \infty.
\]
In the `continuum limit' $k\to \infty$, the discrete sum $\sum_{m\in P\cap k^{-1}\Z^n}$ is replaced by an integral.
Now the RHS of (\ref{FubiniStudypotential}) is to leading order
\[
k^{-1}\log \int_{P} e^{k(-u(m) + \vec{t}\cdot m )} dm, \quad \vec{t}\in \R^n.
\]
This is another Laplace type integral, and its limit as $k\to \infty$ is the Legendre transform of $u$, which gives back the function $\phi$.

The moral is that \emph{in the presence of toric symmetry, algebraic approximation of K\"ahler metrics is related to Legendre transforms}.

\subsection{Extension of K\"ahler currents}

Extension theorems allow us to think extrinsically about K\"ahler currents on subvarieties in some ambient projective manifold.

\begin{thm}\label{extensionKahlercurrent}
(\cite[Thm. B]{Coman}) Let $(X,\omega)$ be a projective manifold with a K\"ahler form representing an integral class, and $Y$ be a smooth subvariety of $X$. Then any $\phi \in PSH(Y,\omega|_Y)$ extends to $\phi\in PSH(X,\omega)$.
\end{thm}

\subsection{Savin's small perturbation theorem}

Savin \cite{Savin} proved that for a large class of second order elliptic equations satisfying certain structural conditions, any viscosity solution $C^0$-close to a given smooth solution has interior $C^{2,\gamma}$-bound. In particular this applies to complex MA equation. Combined with the Schauder estimate,

\begin{thm}\label{Savin}
Fix $k\geq 2$ and  $0<\gamma<1$.
On the unit ball, let $v$ be a given smooth solution to the complex Monge-Amp\`ere equation $(\sqrt{-1}\partial \bar{\partial} v)^n=1$. Then there are constants $0<\epsilon\ll 1$ and $C$ depending on $n, k,\gamma, \norm{v}_{C^{k,\gamma}}$, such that if 
\[
(\sqrt{-1}\partial \bar{\partial} (u+v))^n=1+f, \quad \norm{f}_{C^{k-2,\gamma}}<\epsilon,
\]
and $\norm{u}_{C^0}<\epsilon$, then $\norm{u}_{C^{k,\gamma}(B_{1/2})}\leq C\epsilon$.
\end{thm}



Savin's theorem has fully nonlinear nature, because the perturbative machinery only applies once the solution has \emph{a priori} $C^2$ bound. His proof has two main parts: first he shows a Harnack inequality by a nontrivial application of  Aleksandrov-Bakelman-Pucci estimates, and then uses a compactness argument to prove $C^{2,\gamma}$ estimate, similar to De Giorgi's almost flatness theorem for minimal surfaces.

\subsection{Regularity theory for real Monge-Amp\`ere}\label{RegularitytheoryforrealMA}

There is an extensive literature on the local regularity theory for the real Monge-Amp\`ere equation, largely due to the Caffarelli school. The author thanks C. Mooney for bringing some of these results to his attention. All results surveyed here can be found in \cite{Mooney}.

Any convex function on an open set $v: \Omega\subset \R^n\to \R$ has an associated Borel measure called the Monge-Amp\`ere measure, defined by
\[
MA(v)(E)= |\partial v (E)|,
\]
where $|\partial v(E)|$ denotes the Lebesgue measure of the image of the subgradient map on $E\subset \Omega$. Given a Borel measure $\mu$, a solution to $MA(v)=\mu$ is called an \emph{Aleksandrov solution} to 
$
\det(D^2 v)= \mu;
$ 
if $v\in C^2$, this is the classical real Monge-Amp\`ere equation. 
We shall assume a two-sided density bound
\[
\det(D^2 v)=f \text{ in }B_1, \quad 0<\Lambda_1\leq f\leq \Lambda_2.
\]
Let $B_1\setminus \Sigma$ be the set of strictly convex points of $v$, namely there is a supporting hyperplane touching the graph of $v$ only at one point. Then Caffarelli \cite{Caff1}\cite{Caff2}\cite{Caff4} shows
\begin{itemize}
\item
If $f\in C^\gamma(B_1)$, then $v\in C^{2,\gamma}_{loc}(B_1\setminus \Sigma)$. Then by Schauder theory, if $f$ is smooth, then $v$ is smooth in $B_1\setminus \Sigma$.
\item
If $L$ is a supporting affine linear function to $v$, such that the convex set $\{v=L\}$ is not a point. Then $\{v=L\}$ has no extremal point in the interior of $B_1$. 
\item
The above affine linear set $\{v=L\}$ has dimension $k< n/2$. 
\end{itemize}

Mooney \cite{Mooney} shows further that
\begin{itemize}
\item The singular set $\Sigma$ has $(n-1)$-Hausdorff measure zero. Consequently $B_1\setminus \Sigma$ is path connected (because a generic path joining two given points does not intersect a subset of zero $(n-1)$-Hausdorff measure).

\item  The solution $v\in W^{2,1}_{loc}(B_1)$ even if $\Sigma$ is nonempty.
\end{itemize}

\begin{rmk}\label{Mooneyremark}
A classical counterexample of Pogorelov shows that for $n=3$, the singular set $\Sigma$ can contain a line segment. This is generalised by Caffarelli \cite{Caff4}, who for any $k<n/2$ constructs examples where $f$ is smooth but $\Sigma$ contains a $k$-plane. A surprising example of Mooney \cite{Mooney} shows that the Hausdorff dimension of $\Sigma$ can be larger than $n-1-\epsilon$ for any small $\epsilon$. This means the local regularity theory surveyed above is essentially optimal.
\end{rmk}

\begin{rmk}
On a compact Hessian manifold, the real MA equation makes sense, and Viaclovsky and Caffarelli \cite{CaffarelliViaclovsky} show that the interior singularity cannot occur if the density $f$ is smooth and positive. 
\end{rmk}

\subsection{Special Lagrangian fibration}\label{SpecialLagrangiansurvey}

A real $n$-dimensional submanifold $L$ of a compact Calabi-Yau n-fold $(X, \omega, J, \Omega)$ is called a \emph{special Lagrangian} (SLag) with phase angle $\theta$ if
\begin{equation}
\omega|_L=0,\quad \text{Im}(e^{i\theta}\Omega)|_L=0.
\end{equation}
They are special cases of calibrated submanifolds introduced by Harvey and Lawson \cite{HarveyLawson}, and in particular are minimal submanifolds.
The classical result of McLean says that the deformation theory of SLags with phase $\theta$ is unobstructed, and the first order deformation space is isomorphic to $H^1(L, \R)$. Thus if $L$ is diffeomorphic to $T^n$, then the deformation space is $n$-dimensional, compatible with the SYZ conjecture that $X$ admits a SLag $T^n$-fibration. A sufficient condition to construct Slag fibrations, under the very strong hypothesis of collapsing metric with locally bounded sectional curvature, is obtained by Zhang \cite[Thm 1.1]{Zhang}.

The essence of Zhang's result is a standard application of the implicit function theorem, and we shall summarize the key points (\cf \cite[section 4]{Zhang} for more details). Denote $Y_r= T^n\times B(0, r)\subset T^n_{x_i}\times\R^n_{y_i}\simeq T^* T^n$, where $r\gg 1$ is fixed. The trivial example of a SLag fibration is the following: the CY structure is the flat model
\[
g= \sum (dx_i^2+ dy_i^2), \quad \omega= \sum dx_i \wedge dy_i, \quad \Omega= \bigwedge (dx_j+ \sqrt{-1}dy_j),
\]
and the Slag fibration is just the projection to the $\R^n_{y_i}$ factor, namely the tori $T^n\times \{y\}$ are SLags. Zhang considers a family of CY structures $(g_k, \omega_k, \Omega_k)$ converging to $(g, \omega, \Omega)$ in the $C^\infty$-sense on $Y_{2r}$ (which follows from his bounded sectional curvature assumptions by elliptic bootstrap), such that $\omega_k\in [\omega]\in H^2(Y_{2r}, \R)$. Small deformations of the standard $T^n$ fibres can be represented as graphs on $T^n$: for $y\in \R^n$ and a 1-form $\sigma$ on $T^n$ orthogonal to the harmonic 1-forms $dx_1, \ldots dx_n$, write
\[
L(y, \sigma)= \text{Graph}( x\mapsto y+ \sigma(x)  ) \subset T^*T^n.
\]
The condition for $L(y, \sigma)$ to be a SLag with respect to $(g_k, \omega_k, \Omega_k)$ is
\begin{equation}\label{SlagZhang}
\omega_k|_{ L(y,\sigma) }=0, \quad \text{Im}( e^{\sqrt{-1}\theta_k} \Omega_k  )|_{ L(y,\sigma)  }=0,
\end{equation}
where $\theta_k$ are chosen so that $\int_{T^n}  e^{\sqrt{-1}\theta_k}\Omega_k >0$.  Zhang shows by perturbation arguments that for each $y\in B(0, \frac{3r}{2})$ and $k\geq k_0\gg 1$, there is a unique $\sigma=\sigma_{k,y}$ such that $L(y, \sigma_{k,y})$ solves (\ref{SlagZhang}) with small norm bound $\norm{\sigma_{k,y}}<\delta\ll 1$. He then uses another implicit function argument to show that these SLags indeed define a local SLag $T^n$-fibration on some open subset of $Y_{3r/2}$ containing $Y_r$.


\section{Degenerating Calabi-Yau hypersurfaces}\label{Degeneratingtorichypersurfaces}

We now set the scene for the main work: a particular class of Calabi-Yau hypersurfaces $X_s$ inside
$\mathbb{CP}^{n+1}$ near the large complex structure limit, polarised by the class $\mathcal{O}(n+2)|_{X_s}$ up to a rescaling factor. Special attention will be focused on the simplest case of the Fermat family (\cf Example \ref{Fermatfamily}). We freely borrow from Haase-Zharkov \cite{HaaseZharkov}\cite{HaaseZharkov2}, whose setting includes more general CY hypersurfaces in toric varieties. The key notion is that the degenerating complex structures are controlled by \emph{piecewise linear} data, an idea  studied extensively under the name of tropical geometry.

The philosophy is that every concept in K\"ahler geometry ought to have an analogue in the tropical world, and the combinatorial nature of the tropical version should simplify the original problem in K\"ahler geometry. However, it does not appear clear what is the tropical analogue of the notion of K\"ahler metrics; we devote section \ref{Extensionpropertylocallyconvexsection} and \ref{Fermatcaseextensionpropertysection} to investigate this question, and answer it in the Fermat case by utilizing the large discrete symmetry group.

\subsection{Complex structure}\label{Complexstructure}

Let $N\simeq \Z^{n+1}$, and $M=\Hom(N,\Z)$, and denote $N_\R=N\otimes \R$, $M_\R=M\otimes \R$.
We regard $\mathbb{CP}^{n+1}$ as a toric Fano manifold $\mathbb{P}_\Delta$, with moment  polytope $\Delta\subset M_\R$ corresponding to the anticanonical class $\mathcal{O}(n+2)$. More explicitly $\Delta$ is the $(n+1)$-simplex inside $ M_\R\simeq \{ \sum_{0}^{n+1}y_i=0 \}\subset \R^{n+2}$ spanned by the vertices
\[
(n+1, -1,\ldots -1), (-1, n+1, -1, \ldots -1), \ldots, (-1,\ldots,-1,n+1);
\]
in particular $\Delta$ is a reflexive integral Delzant polytope, with dual polytope
\[
\Delta^\vee= \{ w\in N\otimes \R| \langle m, w \rangle\geq -1, \forall m\in \Delta  \}   \subset \R^{n+2}/\R(1,1,\ldots 1)     
\]
being the $(n+1)$-simplex spanned by the vertices $(1,0,\ldots,0), \ldots, (0,\ldots,0,1)$.
The integral points 
 $m\in \Delta_\Z=\Delta\cap M$ parametrize monomials $z^m$ in the anticanonical linear system $H^0(\mathbb{P}_\Delta, \mathcal{O}(n+2))$.  We study the family of hypersurfaces 
\begin{equation}\label{hypersurface}
X_s= \{ F_s(z)= \sum_{ m\in \Delta_\Z } a_m e^{s\lambda(m)}   z^m =0           \}\subset \mathbb{P}_\Delta, \quad s \gg 1.
\end{equation}
Here $a_m$ are a fixed collection of  coefficients, with $a_0=1$ corresponding to the unique interior integral point $0\in \Delta_\Z$. For any vertex $m$ of $\Delta$, we require $a_m\neq 0$. The function $\lambda$ is defined for those $m\in \Delta_\Z$ for which $a_m\neq 0$; by assumption $\lambda(0)=0$, and $\lambda(m)<0$ otherwise. The natural piecewise linear extension of $\lambda$ to $M_\R$ is assumed to be concave, whose domains of linearity are by assumption simplices, producing a triangulation of $\Delta$. Using the adjunction formula, we can write down a holomorphic volume form $\Omega_s$ on $X_s$, such that along $X_s$
\begin{equation}\label{holomorphicvolume}
dF_s \wedge \Omega_s= d\log z^1\wedge \ldots d\log z^{n+1} ,
\end{equation}
with $z^1, z^2, \ldots z^{n+1}$ the standard coordinates on $(\C^*)^{n+1}\subset \mathbb{P}_\Delta$.
We will always assume $s\gg 1$, and all  the constants in the estimates are independent of $s$.

\begin{eg}\label{Fermatfamily}
The \emph{Fermat family} is given explicitly as
\begin{equation}
X_s= \{ Z_0Z_1\ldots Z_{n+1}+ e^{-s} \sum_{i=0}^{n+1} Z_i^{n+2} =0       \},
\end{equation}
namely we choose $a_m=1$ for $m$ corresponding to the monomials $ Z_0\ldots Z_{n+1}$ and $ Z_i^{n+2}$, and choose $\lambda$ to be the piecewise linear function with value $0$ at the origin and $-1$ at the vertices of $\Delta$.
\end{eg}


The key notion to describe the complex structure degeneration is a piecewise linear object called the \emph{tropicalisation} of the hypersurfaces. Define the nonnegative piecewise linear function $L_\lambda$ on $N_\R$ by
\[
L_\lambda (x)= \max_{m\in \Delta_\Z, a_m\neq 0  } \{ \langle x, m\rangle  +\lambda(m)   \}.
\]
The tropicalisation $\mathcal{A}_\lambda^\infty$ is defined as the nonsmooth locus of $L_\lambda$, or equivalently the locus inside $N_\R$ where the maximum $L_\lambda$ is achieved by at least two values of $m$. There is precisely one bounded component in the complement of $\mathcal{A}_\lambda^\infty$,
\[
\Delta_\lambda^\vee= \{   x| L_\lambda(x)=0    \} \subset N_\R,
\]
whose boundary $\partial \Delta_\lambda^\vee \subset \mathcal{A}_\lambda^\infty$. The relation between the hypersurfaces and the tropicalisation is furnished by the \emph{rescaled log map},
\[
\text{Log}_s: \mathbb{P}_\Delta \supset (\C^*)^{n+1}\to \R^{n+1}\simeq N_\R, \quad 
\text{Log}_s (z)= \frac{1}{s }(\log |z_1|, \ldots \log |z_{n+1}| ).
\]
The image $\mathcal{A}_\lambda^s=\text{Log}_s  (X_s\cap (\C^*)^{n+1})$ is called the \emph{amoeba}. The following Prop. will be tacitly used frequently, as it allows us to think of regions on $X_s$ efficiently in terms of the regions on $\mathcal{A}_\lambda^\infty$, up to a tiny amount of fuzziness.

\begin{prop}\label{amoebaHausdorffconvergence}
(\cf  \cite[Prop. 3.2]{HaaseZharkov}) The amoebas $\mathcal{A}_\lambda^s$ converge in the Hausdorff distance to $\mathcal{A}_\lambda^\infty$ in the Hausdorff distance as $s\to \infty$. In fact 
\[
\begin{cases}
\text{dist}_{\R^{n+1}} (x, \mathcal{A}_\lambda^\infty) \leq \frac{C}{ s  }, \quad \forall x\in \text{Log}_s(X_s),
\\
\text{dist}_{\R^{n+1}} (x, \text{Log}_s(X_s)  ) \leq \frac{C}{ s  }, \quad \forall x\in  \mathcal{A}_\lambda^\infty.
\end{cases}
\]
\end{prop}

\begin{proof}
(sketch) Let $x= \text{Log}_s(z)$ and let $m'\in \Delta_\Z$ saturate the maximum for $L_\lambda(x)$. Applying $\text{Log}_s$ to the inequality
\[
 |e^{s\lambda(m')}z^{m'}|= |-\sum_{m\neq m'} \frac{ a_m}{a_{m'}} e^{s\lambda(m)} z^m| \leq C \max_{m\neq m'}\{  e^{s\lambda(m)} |z^m|    \},
\]
we see 
\[
L_\lambda(x)= \langle x, m'\rangle +\lambda(m') \leq  \max_{m\neq m'}\{ \langle x, m\rangle +\lambda(m)\} + \frac{C}{s},
\]
so  $\text{dist}_{\R^{n+1}} (x, \mathcal{A}_\lambda^\infty) \leq \frac{C}{ s  }$. The other inequality of the claim can be proved by constructing local models of $X_s$ in regions whose $\text{Log}_s$-images are close to $x\in \mathcal{A}_\lambda^\infty$, and then use the implicit function theorem to show $X_s$ is a small perturbation of these local models. 
\end{proof}

\begin{eg}
In the Fermat family example above $\Delta_\lambda^\vee=- \Delta^\vee$ is the reflexion of $\Delta^\vee$.
\end{eg}

The tropicalisation $\mathcal{A}_\lambda^\infty$ is naturally \emph{stratified} according to the subset of $m\in \Delta_\Z$ saturating the maximum $L_\lambda(x)$. This induces a kind of quantitative stratification structure on $\mathcal{A}_\lambda^s$ for $s\gg 1$.

\begin{lem}\label{maximumonsimplex}
There is a fixed number $\delta_1>0$, such that for $s\gg 1$ and any $x\in \mathcal{A}^s_\lambda$ (or $x\in\mathcal{A}^\infty_\lambda$), there is a simplex $\sigma$ in the triangulation of $\Delta$, verifying $\langle x, m\rangle+ \lambda(m)< L_\lambda(x)-\delta_1$ for $m\in \Delta_\Z\setminus \sigma$. 
\end{lem}

\begin{proof}
(Sketch) For any fixed $x\in N_\R$, the function $\langle x,m\rangle+\lambda(m)$ is a concave function of $m\in M_\R$. By our assumptions, the set of $m\in \Delta_\Z$ saturating the maximum must be the set of vertices of some simplex $\sigma$ in the triangulation of $\Delta$. A more effective version of this observation is the Lemma in the $\mathcal{A}^\infty_\lambda$ case, and the $\mathcal{A}^s_\lambda$ case follows by Prop. \ref{amoebaHausdorffconvergence}.
\end{proof}

Given a simplex $\sigma\subset \partial \Delta$ in the triangulation, we associate a subset $\mathcal{A}_{\lambda,\sigma}^\infty$:
\[
\mathcal{A}_{\lambda,\sigma}^\infty=\{     x\in \mathcal{A}_\lambda^\infty | L_\lambda(x)= \langle x, m\rangle +\lambda(m),  \forall m\in \sigma   \} .
\]
Clearly if $\sigma \prec \sigma'$, then $\mathcal{A}_{\lambda,\sigma}^\infty \supset \mathcal{A}_{\lambda,\sigma'}^\infty$.  The intuition is that larger $\sigma$ correspond to more nongeneric regions, and the complement of their neighbourhoods correspond to more generic regions.

\begin{Notation}
We need a few terminologies to describe $\mathcal{A}_{\lambda,\sigma}^\infty$. The face of $\partial\Delta_\lambda^\vee$ dual to $\sigma$ is $F_\sigma^\vee= \partial\Delta_\lambda^\vee \cap \mathcal{A}_{\lambda,\sigma}^\infty$. The outward normal cone to $\sigma$ is
\[
NC_\Delta(\sigma)= \{  x\in N_\R | \langle m, x\rangle \leq  \langle m', x\rangle, \forall m\in \Delta_\Z, \forall m' \in \sigma    \}.
\]
By the Delzant polytope property $NC_\Delta(\sigma)$ is isomorphic to $\R_{\geq 0}^l$, where $n+1-l$ is the dimension of the minimal face of $\partial \Delta$ containing $\sigma$.
The Minkowski sum of two sets $A,B$ means $A+B=\{ a+b|a\in A, b\in B  \}$. 	
\end{Notation}

\begin{lem}\label{stratasigma}\label{Minkowskisumlemma}
	(compare \cite[Lemma 3.1]{HaaseZharkov}) If $\dim\sigma\geq 1$ then
$
\mathcal{A}_{\lambda,\sigma}^\infty= F_\sigma^\vee+  NC_\Delta(\sigma).
$
\end{lem}

\begin{lem}
$\mathcal{A}_\lambda^\infty=\partial \Delta_\lambda^\vee\cup  \bigcup_{   \dim \sigma \geq 1} \mathcal{A}_{\lambda,\sigma}^\infty  $.
\end{lem}

\begin{proof} Let $x\in \mathcal{A}_\lambda^\infty$.
If $m=0\in \Delta$ achieves the maximum $L_\lambda(x)$, then $x\in \partial \Delta_\lambda^\vee$. If not, then the
 maximum is achieved by at least two $m\in \partial \Delta$, so $x\in \mathcal{A}_{\lambda,\sigma}^\infty$ for some $\sigma \subset \partial \Delta$ with $\dim \sigma\geq 1$. 
\end{proof}

\begin{rmk}
The intuition is that a neighbourhood of $\partial \Delta_\lambda^\vee$ corresponds to a toric region, while $\mathcal{A}_{\lambda,\sigma}^\infty$ controls how $X_s$ approaches the toric boundary of $\mathbb{P}_\Delta$, and the stratification is related to how the toric boundary components intersect.	
\end{rmk}

Our next goal is to assign good holomorphic charts to $X_s$ related to the stratification structure. We first consider the \emph{toric region}, which  shall be covered by $(\C^*)^n$-charts. Let $w\in N$ be the primitive integral outward normal vector to a facet $F(w)=\{ m\in \Delta| \langle w, m\rangle=1  \}$ of $\Delta$. The chart parametrised by $w$ is contained inside the region 
\begin{equation}
U_w^{s, o}= \{ z\in X_s|  e^{s\lambda(m) } |z^m| \ll 1, \quad \forall m\in \Delta_\Z \setminus (F(w)\cup \{0\} )  \}.
\end{equation}
Let $m_0\in F(w)\cap \Delta_\Z$, and choose an integral basis $m_1, \ldots, m_n$ for $\{ m\in M| \langle w, m\rangle =0  \}$. Then the monomials $z^{m_1},\ldots, z^{m_n}$ provide the local $(\C^*)^n$-coordinates on the chart, since by the implicit function theorem $X_s$ is locally a graph $\{  z^{m_0}= f(z^{m_1}, \ldots , z^{m_n})  \}$. In fact by the defining equation (\ref{hypersurface}) of the hypersurface
\[
z^{-m_0} \approx -\sum_{ m \in F(w) }  a_m  e^{s\lambda(m)} z^{m-m_0} ,
\]
whence the holomorphic volume form is (\cf (\ref{holomorphicvolume}))
\begin{equation}\label{holomorphicvolumetoricregion}
\Omega_s= \pm\frac{ d\log z^{m_0} \wedge \ldots d\log z^{m_n} }{dF_s } \approx d\log z^{m_1} \wedge \ldots d\log z^{m_n}.
\end{equation}
(Here $m_i$ are suitably oriented to  take care of $\pm1$.) We regard the above region as an open subset of $(\C^*)^n$, and denote the chart $U^s_w$ as the largest $T^n$-invariant subset, delineated by a collection of affine linear inequalities on the variables $\log |z^{m_i}|$.

In the tropical limit $s=\infty$, the region $\text{Log}_s(U_w^{s,o})$ becomes 
\[
U_w^{\infty, o}= \{    x\in \mathcal{A}_\lambda^\infty | \langle m, x\rangle+ \lambda(m)<0, \forall m\in \Delta_\Z \setminus (F(w)\cup \{0\})           \}
\]
Inside this the limiting version of $\text{Log}_s(U_w^{s})$ is
\[
U_{w}^\infty= \{ (U_w^{\infty,o}\cap \partial \Delta_\lambda^\vee) + \R_{\geq 0} w\} \cap \mathcal{A}_\lambda^\infty  .
\]
Later we shall also need the slightly shrinked regions for $0<\delta\ll \delta_1$: let
\[
U^{s,o}_{w,\delta}=   \{ z\in X_s|  e^{s\lambda(m) } |z^m| \ll e^{-s\delta}, \quad \forall m\in \Delta_\Z \setminus (F(w)\cup \{0\} )  \},
\]
whose largest $T^n$-invariant subset is 
$
U^s_{w, \delta} $. The tropical limit of $\text{Log}_s(U_{w,\delta}^{s,o})$ is
\[
U^{\infty,o}_{w,\delta}= \{    x\in \mathcal{A}_\lambda^\infty| \langle m, x\rangle+ \lambda(m)<-\delta, \forall m\in \Delta_\Z \setminus (F(w)\cup \{0\})           \},
\]
containing the limiting version of $\text{Log}_s(U_{w,\delta}^{s})$
\[
U^\infty_{w,\delta}= \{ (U_w^{\infty,o}\cap \partial \Delta_\lambda^\vee) + \R_{\geq 0} w\} \cap \mathcal{A}_\lambda^\infty  .
\]
As the choice of $w$ varies, such regions $U_{w,\delta}^{\infty,o}$ cover a neighbourhood of $\partial \Delta_\lambda^\vee$ as a consequence of Lemma \ref{maximumonsimplex}; so do $U_{w,\delta}^{\infty}$. This means the charts of toric type already cover part of the neighbourhood of the toric boundary.

\begin{eg}
In the $n=1$ case, $X_s$ are elliptic curves, and the toric charts cover the entire $X_s$. In the $n=2$ case, $X_s$ are quartic K3 surfaces, and the toric charts cover most parts of $X_s$ including a large portion of the intersection of $X_s$ with the toric boundary of $\mathbb{P}^3$, but do not cover a tiny neighbourhood of the 24 points located at the intersection of $X_s$ with $\{  Z_i=Z_j=0 \}$.
\end{eg}

We now consider the neighbourhood of the \emph{toric boundary} near the stratum $\mathcal{A}_{\lambda,\sigma}^\infty$, but keeping away from higher strata and from $\partial \Delta_\lambda^\vee$. 
Here
\[
|a_{m'} e^{s\lambda(m')} z^{m'}|\ll |a_m e^{s\lambda(m)} z^m|, \quad \forall m'\in \Delta_\Z\setminus \sigma, \forall m\in \sigma.
\]
 Since 
 most terms in the defining equation (\ref{hypersurface}) are negligible in our region, the hypersurface is locally approximately
 \[
 \sum_{ m\in \sigma\cap \Delta_\Z} a_m e^{s\lambda(m)} z^m \approx 0.
 \]
We focus on the subregion where $m_0'\in \sigma$ achieves the maximal magnitude for $|a_m e^{s\lambda(m)} z^m|$, and $m_1'\in \sigma$ achieves the second largest magnitude. These two magnitudes must be of comparable size by the hypersurface equation. Choose an integral basis $w_1,\ldots w_l$ for the outward normal cone $NC_\Delta(\sigma)$,
so $\langle m, w_i \rangle =1$ for $m\in \sigma\cap \Delta_\Z$. Denote the vertices of $\sigma$ as $m_i'$ for $i=0,1,\ldots, \dim \sigma$, and choose $m_1, \ldots m_{\dim \sigma-1}$ an integral basis of $\text{span}_\Q\{ m_2'-m_0', \ldots, m_{\dim \sigma}'-m_0'\}\cap M$. Choose $m_0$ so that $m_0, \ldots m_{\dim \sigma-1}$ is an integral basis of $\text{span}_\Q\{ m_1'-m_0', \ldots, m_{\dim \sigma}'-m_0'\}\cap M$, and complete this into an integral basis $\{ m_0, \ldots, m_{n-l}  \}$ for $\text{span}\{w_1,\ldots w_l\}^\perp $, providing $(n+1-l)$ $\C^*$-variables $z^{m_0}, \ldots z^{m_{n-l}}$. We then find $\mathfrak{m}_j$ for $j=1,2,\ldots l$, with $\langle \mathfrak{m}_j, w_i\rangle =-\delta_{ij}$, and we can demand $\mathfrak{m}_1+\ldots \mathfrak{m}_l=-m_0'$ because $\langle m_0', w_i\rangle= 1$. These provide the $\C$-variables $z^{\mathfrak{m}_j}$ for $1\leq j\leq l$, which can vanish on the toric boundary. On this local piece of $X_s$, 
the variables $z^{\mathfrak{m}_j}$ and $z^{m_1},\ldots z^{m_{n-l}}$ furnish a set of local coordinates as the $\C^*$-variable $z^{m_0}$ is expressible locally as a function of theirs.

The holomorphic volume form (\ref{holomorphicvolume}) is 
\begin{equation}\label{holomorphicvolumetoricboundary}
\begin{split}
\Omega_s= & \pm \frac{ d\log z^{ \mathfrak{m}_1 } \wedge \ldots d\log z^{\mathfrak{m}_l}\wedge d\log z^{m_0}\wedge \ldots d\log z^{m_{n-l}} }{dF_s }
\\
=& \pm
\frac{dz^{\mathfrak{m}_1}\wedge \ldots dz^{\mathfrak{m}_l} }{ z^{-m_0'}  } \bigwedge    \frac{ d\log z^{m_0}\wedge \ldots d\log z^{m_{n-l}} }{dF_s } 
\\
\approx & \pm
dz^{\mathfrak{m}_1}\wedge \ldots dz^{\mathfrak{m}_l}  \bigwedge    d\log z^{m_1}\wedge \ldots d\log z^{m_{n-l}}  \bigwedge \frac{d\log z^{m_0} }{ a_{m_1'} e^{s\lambda(m_1') } dz^{m_1'-m_0'} }
\\
= &  \frac{dz^{\mathfrak{m}_1}\wedge \ldots dz^{\mathfrak{m}_l} }{ a_{m_1'} e^{s\lambda(m_1')  } z^{m_1'-m_0'} \mathfrak{d} } \bigwedge d\log z^{m_1} \wedge \ldots d\log z^{m_{n-l}},
\end{split}
\end{equation}
up to choosing appropriate ordering of the coordinates. Here $\mathfrak{d}$ is the divisibility of $m_1'-m_0'$ inside the group
\[
\text{span}_\Q\{  m_1'-m_0', \ldots, m_{\dim \sigma}'-m_0'    \}\cap M/ \text{span}_\Z\{  m_1, \ldots, m_{\dim \sigma-1}    \} \simeq \Z.
\]
Notice $a_{m_1'} e^{s\lambda(m_1')  } z^{m_1'-m_0'}$ is uniformly equivalent to $a_{m_0'} e^{s\lambda(m_0')}$ in this region.

\begin{rmk}
The discussion above can be simplified if one assumes the triangulation of $\Delta$ is \emph{maximal}, namely each simplex is $\Z$-isomorphic to the standard simplex. We choose not to do so because this stronger assumption would exclude the Fermat family.
\end{rmk}

\begin{rmk}\label{chartsofboundarytype}
A problem when we work with the coordinates $z^{m_1}, \ldots z^{m_{n-l}}, z^{\mathfrak{m}_j}$ is the inequality constraint to keep $a_{m_0'} e^{s\lambda(m_0')} z^{m_0'}$ and $a_{m_1'} e^{s\lambda(m_1')} z^{m_1'}$ as  the two dominant monomials. This means such a holomorphic chart is not quite as simple as the product of $D(1)^\ell$ with a long annulus in $(\C^*)^{n-l}$. In practice we will cover this region by lots of simpler charts which we call the \emph{charts of boundary type}. Let $P$ be any point in this region, such that $\max_m{ |a_m  e^{s\lambda (m)} z^m|  }$ is large but still comparable to 1 (to guarantee the chart overlaps nontrivially with some toric type chart). The associated chart uses the same coordinates $z^{m_1}, \ldots z^{m_{n-l}}, z^{\mathfrak{m}_j}$ as above, but describes only a small region:
\[
U_P= \{   |   z^{\mathfrak{m}_j }| \lesssim  |   z^{\mathfrak{m}_j}(P) |, \forall j, \quad 
|z^{m_i}- z^{m_i}(P)|< c |z^{m_i}(P)|,   \forall i   \} ,
\]
where $0<c\ll 1$ is a fixed dimensional constant. These charts have an interpretation in terms of the strata $\mathcal{A}_{\lambda,\sigma}^\infty$ (\cf Lemma \ref{Minkowskisumlemma}): the point $P$ corresponds roughly to a point $P'$ on the face $F_\sigma^\vee\subset \Delta_\lambda^\vee$, and allowing $|   z^{\mathfrak{m}_j }|$ to decrease to zero corresponds to taking the Minkowski sum with the outward normal cone $NC_\Delta(\sigma)$, so the tropical analogue of our small chart is $\{ P'\}+ NC_\Delta(\Sigma)$.
 
\end{rmk}

\begin{eg}
	For generic quartic K3 surfaces, the following simple situation models a small neighbourhood of the 24 points on $\text{K3}\cap \{Z_i=Z_j=0\}$. Locally the dominant monomials are $(z_1z_2)^{-1}, (z_1z_2)^{-1}z_0, 1$, where $z_1, z_2$ are $\C$-coodinates which vanish on toric boundaries, and $z_0$ is a $\C^*$-coordinate; together $z_0, z_1, z_2$ are local coordinates on $\mathbb{P}^3$. The local model hypersurface is
	\[
	\{ -(z_1z_2)^{-1} + (z_1z_2)^{-1} z_0 = 1 \}=   \{  z_0 =1+ z_1z_2       \} ,
	\]
	so $z_1, z_2$ can be used as local coordinates on the hypersurface.
	The holomorphic volume form $\Omega$ on the hypersurface is (up to a normalising factor)
	\[
	\Omega= \frac{ d\log z_0 \wedge d\log z_1\wedge d\log z_2}{ d(  -(z_1z_2)^{-1} + (z_1z_2)^{-1} z_0 - 1   )  }=  z_0^{-1} d z_1\wedge d z_2.
	\]
	This is the typical boundary type behaviour. A significant part of the boundary type region overlaps with the toric region. In this example, when $|z_1|$ is not too small, we can view $z_1$ as a $\C^*$-coordinates, so $\{ z_1,z_0 \}$ provides a toric type chart, as we can express $z_2= z_1^{-1}(z_0-1)$. In this chart 
	\[
	\Omega= z_0^{-1}dz_1\wedge dz_2= d\log z_1\wedge d\log z_0,
	\] 
	which agrees with the standard holomorphic volume form in toric type charts. The same behaviour happens when $|z_2|$ is not too small. The problem mentioned in Remark \ref{chartsofboundarytype} is due to the fact that this local model is only a valid approximate description of the K3 for $z_0, z_1, z_2$ satisfying some inequality constraints. The prescription of charts of boundary type means that we are simultaneously using the charts $\{ |z_1|\lesssim \nu, |z_2|\lesssim \nu^{-1}   \}$ for many choices of parameters $\nu$. Notice the scaling symmetry \[
	z_1\mapsto \nu z_1, \quad z_2\mapsto \nu^{-1} z_2
	\]
	means that there is no obviously preferred chart of boundary type. More concrete examples can be found in \cite[section 1.1.6]{Li}.
\end{eg}

Local charts of the toric type and the boundary type cover the entire hypersurface $X_s$ for $s\gg 1$, and a substantial portion of any boundary type chart is in fact already covered by toric charts. Almost all the measure is contained in the toric type region.



\subsection{Piecewise linear structure}\label{Piecewiseaffinestructure}

\begin{prop}
	The polyhedral complex $\partial\Delta_\lambda^\vee$ is homeomorphic to $S^n$.
\end{prop}

\begin{proof}
	This is because $\partial \Delta_\lambda^\vee$ is the boundary of a convex polyhedron $\Delta_\lambda^\vee$ with nontrivial interior.
\end{proof}

We now assign a a collection of charts 
to $\partial \Delta_\lambda^\vee$,  whose transition functions are piecewise linear. (Some authors prefer the terminology `piecewise affine'.) These are closely related to the holomorphic charts on $X_s$ in section \ref{Complexstructure}.

Let $w\in N$ be the primitive integral outward normal vector to a facet $F(w)$ of $\Delta$, and choose an integral basis $m_1,\ldots m_n$ for $\{ m\in M|\langle w, m\rangle=0  \}$, suitably oriented to be compatible with (\ref{holomorphicvolumetoricregion}).  On the open subset of $\partial \Delta_\lambda^\vee$,
\[
\partial\Delta_\lambda^\vee\cap U_w^\infty=
\{ x\in  \partial\Delta_\lambda^\vee  |   \langle m, x\rangle +\lambda(m) <0 , \forall m\in \Delta_\Z\setminus (F(w)\cup \{0\}  )    \},
\]
we regard $m_1, \ldots m_n$ as the affine linear coordinates, also written as $x^{m_1},\ldots x^{m_n}$. Such charts cover $\partial \Delta_\lambda^\vee$. We denote $\widetilde{Sing}$ as the subset of points on $\partial \Delta_\lambda^\vee$ which do not lie on the interior of the top dimensional faces. It is easy to check that the transition functions on overlapping charts in $\partial \Delta_\lambda^\vee\setminus \widetilde{Sing}$ lie in $SL(n,\Z)\ltimes \R^n$, so the volume form $dx^{m_1} \wedge \ldots dx^{m_n}$ is defined independent of the choice of charts. We call the associated measure $d\mu_\infty$ the \emph{Lebesgue measure} on $\partial \Delta_\lambda^\vee$, with respect to which $\widetilde{Sing}$ is a null set. The set $\widetilde{Sing}$ has real codimension 1, and the transition functions are in general only piecewise linear.

\begin{rmk}
The affine structure on $\partial \Delta_\lambda^\vee\setminus \widetilde{Sing}$ can be often extended to a subset of $\partial \Delta_\lambda^\vee$ with codimension 2 complement. This in general involves a somewhat ad hoc choice of the singular locus. In the Fermat family case, due to the discrete symmetry, the barycentric subdivision provides a canonical choice. (\cf section \ref{Fermatcaseextensionpropertysection}).

\end{rmk}



We now examine the \emph{normalised canonical measure} on $X_s$
\begin{equation}\label{normalisedcanonicalmeasure}
d\mu_s =\frac{1}{(4\pi s)^n} \sqrt{-1}^{n^2} \Omega_s\wedge \overline{\Omega}_s.
\end{equation}

\begin{prop}\label{canonicalmeasureconvergence}
	As $s\to +\infty$, the pushforward measure $(\text{Log}_s)_* d\mu_s $ converges to the Lebesgue measure $d\mu_\infty$ supported on $\partial\Delta_\lambda^\vee$. In particular 
	\begin{equation}
	\int_{X_s} d\mu_s\to \text{Vol}(\partial\Delta_\lambda^\vee)= \int_{\partial \Delta_\lambda^\vee} d\mu_\infty .
	\end{equation}
Morever, there is a uniform exponential measure decay estimate
	\begin{equation}\label{exponentialmeasuredecayawayfromtoric}
	d\mu_s (   \{ z\in X_s:  \text{dist}_{\R^{n+1}}(\text{Log}_s (z), \partial \Delta_\lambda^\vee ) > s^{-1}\Lambda        \}         ) \leq C'e^{- C \Lambda  } , \quad \forall \Lambda>0.
	\end{equation}
\end{prop}

\begin{proof}
(Sketch) Using Lemma \ref{stratasigma} and the holomorphic volume form formula (\ref{holomorphicvolumetoricboundary}), the neighbourhood of the toric boundary near $\mathcal{A}_{\lambda,\sigma}^\infty$ only contributes $O(s^{-l})$ to the normalised measure, where $l=\dim NC_\Delta(\sigma)$. The same lemmas imply (\ref{exponentialmeasuredecayawayfromtoric}) by summing over contributions from boundary type regions. In the toric region corresponding to the neighbourhood of $\partial \Delta_\lambda^\vee$, the convergence of the normalised volume measure follows from Prop. \ref{amoebaHausdorffconvergence} and  formula (\ref{holomorphicvolumetoricregion}).
\end{proof}

\begin{rmk}
The measure convergence holds for much more general degenerating families by the work of Boucksom et al. \cite{Boucksom1}. The fact that the measure is concentrated along $\partial \Delta_\lambda^\vee$ justifies why we focus on $\partial \Delta_\lambda^\vee$ rather than $\mathcal{A}_\lambda^\infty$.
\end{rmk}

\subsection{K\"ahlerian polarisation}\label{Kahlerianpolarisation}

We specify a polarisation class $[\Delta]$ on the toric manifold $\mathbb{CP}^{n+1}=\mathbb{P}_\Delta$. A standard background K\"ahler metric is (a suitable multiple of) the Fubini-Study metric:
\[
\begin{split}
\omega_{FS}&= \frac{\sqrt{-1} (n+2)}{2}\partial \bar{\partial} \log (|Z_0|^2+ \ldots |Z_{n+1}|^2)
\\
&= \frac{\sqrt{-1} (n+2)}{2}\partial \bar{\partial} \log (\sum_{m\in vertex(\Delta) } e^{ \frac{2}{n+2}\langle m, \text{Log}(z) \rangle }  )  .
\end{split}
\]
Our normalisation guarantees that the potential has the asymptotic behaviour
\[
\sup_{z} |\frac{ (n+2)}{2} \log (\sum_{m\in vertex(\Delta) } e^{ \frac{2}{n+2}\langle m, \text{Log}(z) \rangle }  ) -\max_{m\in \Delta} \langle m, \text{Log}(z) \rangle |<+\infty.
\]

A general (singular) K\"ahler metric $\omega_u$ on $(\mathbb{P}_\Delta, [\Delta])$ is given by a relative potential $u\in PSH(X,\omega_{FS})$. Alternatively, one thinks of $\omega_u$ as a collection of local absolute potentials: 
\begin{equation}
\begin{cases}
u_0= u+  \frac{ (n+2)}{2} \log (\sum_{m'\in vertex(\Delta) } e^{ \frac{2}{n+2} \langle m', \text{Log}(z) \rangle }  )  ,
\\
u_{m}= u+  \frac{ (n+2)}{2} \log (\sum_{m'\in vertex(\Delta) } e^{ \frac{2}{n+2}\langle m', \text{Log}(z) \rangle }  ) - \langle m, \text{Log}(z) \rangle  ,
\end{cases}
\end{equation}
where $u_0$ is a local potential in a compact region, and $u_m$ give the local potentials near the toric boundary.

 We call a convex function $u$ on $N_\R=\R^{n+1}$ \emph{admissible} if it satisfies the asymptotic growth condition
 \begin{equation}\label{admissible}
 \sup_x |u(x)- \max_{m\in \Delta} \langle m, x\rangle |< +\infty,
 \end{equation}
 which captures the information of the K\"ahler class.

 \begin{prop}\label{admissibleimpliesextension}
 	A convex function $u$ is admissible if and only if the K\"ahler current defined by the psh function $u\circ \text{Log}$ on $(\C^*)^{n+1}$ extends to a torus invariant K\"ahler current on $(\mathbb{P}_\Delta, [\Delta])$ with continuous local potentials.
 \end{prop}

 \begin{proof}
 	(Sketch) Convex functions on $(\C^*)^n$ correspond to torus invariant psh functions via the log map (\cf Lemma \ref{pshconvexity} below). If $u$ is admissible, then near the toric boundary the appropriate local potential $u_m$ extends continuously over the boundary piece by the growth asymptote assumption and convexity, and the extension remains psh. Conversely, the asymptotic condition is dictated by the local boundedness of $u_m$ near the toric boundary pieces.
 \end{proof}

A general (singular) K\"ahler metric $\omega_\varphi$ on $X_s$ in the polarisation class $s^{-1}[\Delta]$ is given by a  potential $\varphi\in PSH(X_s ,s^{-1}\omega_{FS})$. The normalising factor $s^{-1}$ is aimed at extracting nontrivial limits as $s\to \infty$. We can completely analogous define the local potentials:
\begin{equation}\label{localpotentialscomplexsetting}
\begin{cases}
\varphi_0= \varphi +  \frac{ (n+2)}{2s} \log (\sum_{m'\in vertex{(\Delta)} } e^{ \frac{2}{n+2} \langle m', \text{Log}(z) \rangle }  )  ,
\\
\varphi_{m}= \varphi+  \frac{ (n+2)}{2s} \log (\sum_{m'\in   vertex{(\Delta)}    } e^{ \frac{2}{n+2}\langle m', \text{Log}(z) \rangle }  ) - \langle m, \text{Log}_s(z) \rangle  ,
\end{cases}
\end{equation}
which are by definition psh on respective regions.

In particular, we can represent the Calabi-Yau metric $\omega_{CY,s}$ on $X_s$ by a potential $\varphi_{CY,s}$. The Calabi-Yau condition is
\begin{equation}\label{CalabiYaucondition}
\omega_{CY,s}^n=  a_s s^{-n} d\mu_s,
\end{equation}
where the normalising constant 
\begin{equation}\label{ainfty}
a_s= \frac{ \int_{X_s} [\Delta]^n} { \int_{X_s} d\mu_s   }\to a_\infty = \frac{ \int_{X_s} [\Delta]^n} {  \text{Vol}(\partial\Delta_\lambda^\vee)   } 
\end{equation}
as $s\to +\infty$ (\cf Prop. \ref{canonicalmeasureconvergence}).

\subsection{Extension property and locally convex functions}\label{Extensionpropertylocallyconvexsection}

We now discuss the issue of finding a tropical notion analogous to K\"ahler metrics. The concept of a K\"ahler metric is formulated in terms of a collection of local psh functions $\phi_j$ on overlapping complex charts, whose  differences $\{\phi_i-\phi_j  \}$ represent a given cocycle of local pluriharmonic function. Intuitively, the analogue should be a collection of local convex functions $u_j$ whose differences $\{ u_i-u_j \}$ represent a given cocycle of local affine functions.

To the author's awareness there is no definitive formulation of local convexity on polyhedral sets. In the case of interest, we need to define a class of `locally convex functions' on  $\partial\Delta_\lambda^\vee$. The problem is that on $\widetilde{Sing}\subset \partial\Delta_\lambda^\vee$, the transition functions between different charts are only piecewise linear, so convexity is not invariantly defined. This problem also prevents us from setting up a general global notion of real MA equation on $ \partial\Delta_\lambda^\vee$, which is an essential ingredient in the SYZ conjecture in general. We will attempt to give a special definition in the Fermat case (\cf section \ref{Fermatcaseextensionpropertysection}).

However, the extension theorem \ref{extensionKahlercurrent} provides an alternative viewpoint: (1,1)-type K\"ahler currents can be defined extrinsically. By analogy, we propose that the correct notion should be equivalent to the following

\begin{Def}\label{extensionpropertydef}
A continuous function $u$ on $\partial \Delta_\lambda^\vee$ satisfies the \emph{extension property} if it extends to an admissible convex function on $N_\R=\R^{n+1}$ defined in section \ref{Kahlerianpolarisation}. 
\end{Def}

\begin{eg}
The zero function extends to $L_\lambda$, which is admissible and convex.
\end{eg}

The problem is to make this definition both intrinsic to $\partial \Delta_\lambda^\vee$, and local in nature. We do not fully succeed but shall make some partial progress.

\begin{prop}\label{extensionpropertycharacterisation}
A continuous function $u$ on $\partial \Delta_\lambda^\vee$ satisfies the extension property if and only if for every $x\in \partial\Delta_\lambda^\vee$, there exists $p\in \Delta$, such that for any $y\in \partial\Delta_\lambda^\vee$,
\[
u(y) \geq u(x) + \langle p, y-x\rangle.
\]
\end{prop}

\begin{proof}
The if direction is because the asymptotic growth condition (\ref{admissible}) implies the gradient of $u$ must be contained in $\Delta$.

For the only if direction, we apply the Legendre transform:
\[
u^*(p) = \sup_{x\in  \partial \Delta_\lambda^\vee } \{ \langle x, p\rangle - u(x) \}, \quad p\in \Delta,
\]
and consider a version of the double Legendre transform
\[
u^{**}(x)= \sup_{p\in \Delta} \{ \langle x, p\rangle - u^*(p) \} .
\]
Clearly $u^{**}$ is convex, and admissible by the boundedness of $u^*$, and $u^{**}(x)\leq u(x)$ on $\partial \Delta_\lambda^\vee$ because
\[
\langle x,p\rangle -u^*(p) \leq u(x), \quad \forall p\in \Delta.
\]
Our characterisation precisely ensures $u^{**}(x) \geq u(x)$ on $\partial \Delta_\lambda^\vee$. Then $u^{**}$ provides the canonical extension.
\end{proof}

\begin{rmk}
The above characterisation is not completely intrinsic because it uses the extrinsic pairing $\langle, \rangle: M_\R\times N_\R\to \R$. On the positive side it uses only the value of $u$ on $\partial \Delta_\lambda^\vee$.
\end{rmk}

\begin{rmk}\label{pisgradient}
At $x\in \partial \Delta_\lambda^\vee$, the vector $p\in \Delta$ in the hypothesis is a subgradient of the canonical extension $u^{**}$, namely $u^{**}(y)- u(x)\geq \langle p, y-x\rangle$. 
\end{rmk}

We now introduce a local notion. The function $u$ below will be analogous to $\phi_0$ in (\ref{localpotentialscomplexsetting}). Recall the charts $\partial \Delta_\lambda^\vee\cap U_w^\infty$ associated to ourward normal vectors $w$ introduced in section \ref{Piecewiseaffinestructure}, with local coordinates $x^{m_1}, \ldots x^{m_n}$.

\begin{Def}\label{locallyconvexdef}
Let $u$ be a continuous function on $\partial \Delta^\vee_\lambda$, to which we associate a collection of local functions $\{ u_m \}_{ m\in \Delta_\Z }$ by the rule $u_m=u- \langle x, m\rangle$. We regard $u_m$ as a function on the charts $\partial \Delta_\lambda^\vee\cap U_w^\infty$ with $\langle w, m\rangle = 1$. We say $u$ is a \emph{locally convex} function if all $u_m$ are convex on their corresponding charts. 
\end{Def}

\begin{rmk}
One can reconstruct $u$ from the local functions $\{ u_m \}$ as long as their mutural differences define a correct cocycle $\{m-m' \}$. Thus this definition has the intrinsic local feature we desire, in analogy with the notion of K\"ahler potentials.
\end{rmk}

\begin{rmk}\label{caveatconvexity}
For fixed $w$ and $m,m'$ satisfying $\langle m,w\rangle=\langle m', w\rangle=1$, the convexity of $u_m$ and $u_{m'}$ on the $w$-chart are equivalent because $m-m'$ is an affine function. However, on the overlap of the $w$-chart and the $w'$-chart, if $u_m$ is convex in one chart it is not automatically convex in the other.
\end{rmk}

\begin{rmk}
If a convex function is not sufficiently regular, there can be a null set of points at which the subgradient is not unique. Later we will abuse language to use the word gradient to refer to any choice of subgradient.
\end{rmk}

\begin{prop}\label{extensionimplieslocalconvexity}
If $u$ satisfies the extension property, then $u$ is locally convex.
\end{prop}

\begin{proof}
Let $\langle m, w \rangle=1$, and consider the function $u_m$ on the chart $\partial \Delta_\lambda^\vee\cap U_w^\infty$. Given $x$ in the chart, we need to find $\vec{p}$ such that 
\[
u_m(y)- u_m(x)  \geq \vec{p}\cdot (y-x)_w,
\]
where $\vec{p}$ is a covector, and $(y-x)_w$ refers to the representation of $y-x$ in the local coordinates $x^{m_1},\ldots x^{m_n}$; after identifying $x^{m_1},\ldots x^{m_n}$ as coordinates on the plane $m^\perp= \{ \langle m, x'\rangle=0 \}$, we may regard $(y-x)_w$ as an element of $m^\perp$, and according to the decomposition $N_\R= m^\perp \oplus \R w   $,
\[
y-x= (y-x)_w + \langle y-x, m \rangle w.
\] 
Since the convexity of $u_m$ and $u_{m'}$ are equivalent in the $w$-chart if $\langle m, w\rangle=\langle m, w\rangle=1$, we may assume  $L_\lambda(x)$ is attained by $\langle m, x\rangle+ \lambda(m)$. By the extension property and Prop. \ref{extensionpropertycharacterisation}, there is some $p\in \Delta$, such that 
\[
u(y)- u(x) \geq \langle p, y-x\rangle,
\]
hence
\[
u_m(y)- u_m(x) \geq \langle p-m, y-x\rangle = \langle p-m, (y-x)_w \rangle + \langle y-x, m \rangle \langle p-m, w\rangle.
\]
Since $p\in \Delta$,
we have $
\langle p, w\rangle \leq 1 = \langle m, w\rangle.
$
Since $L_\lambda(x)$ is attained by $\langle m, x\rangle+ \lambda(m)$, and  the polytope $\Delta_\lambda^\vee$ lies in the half space $\{ \langle m, \rangle+ \lambda(m) \leq 0 \}$, we have
\[
\langle m, y\rangle+ \lambda(m) \leq 0= \langle m, x\rangle+ \lambda(m) .
\]
Combining the above
\[
u_m(x)- u_m(y) \geq  \langle p-m, (y-x)_w \rangle + \langle y-x, m \rangle \langle p-m, w\rangle \geq \langle p-m, (y-x)_w \rangle,
\]
so we have produced $\vec{p}$ as required.
\end{proof}

\subsection{Extension property: the Fermat case}\label{Fermatcaseextensionpropertysection}

We do not know the equivalence between the extension property and the local convexity property. However, in the case of the Fermat family Example \ref{Fermatfamily}, the polyhedral set $\partial \Delta_\lambda^\vee=- \partial \Delta^\vee$ has a \emph{discrete symmetry} by the permutation group of the vertices of $\Delta$, corresponding to the permutations of the monomials $Z_0^{n+2}, \ldots Z_{n+1}^{n+2}$. This can be used to our advantage.

\begin{Notation}
Denote the vertices of $\partial \Delta_\lambda^\vee$ as $w_0, \ldots, w_{n+1}$, which coincide with the outward normal vectors because $\partial \Delta_\lambda^\vee=- \partial \Delta^\vee$. Denote the vertices of $\Delta$ as $m^0, \ldots, m^{n+1}$, so that 
\[
\langle w_i, m^j\rangle =\begin{cases}
1, \quad & i\neq j, \\
-(n+1), \quad & i=j.
\end{cases}
\]
Let $\text{Star}(w_i)$ be the star of $w_i$ in the barycentric subdivision of $\partial \Delta_\lambda^\vee$. Let $Sing \subset \widetilde{Sing}$ be the subset of points not contained in the interior of any of these stars. The affine structure on $\partial \Delta_\lambda^\vee\setminus \widetilde{Sing}$ extends to $\partial \Delta_\lambda^\vee\setminus Sing$, by decreeing that on the interior of $\text{Star}(w_i)$ we use the coordinates for the chart $U_w^\infty\cap \partial \Delta_\lambda^\vee$.  As $Sing$ has codimension two inside $\partial \Delta_\lambda^\vee$, this makes $\partial \Delta_\lambda^\vee$ into a \emph{singular affine manifold}.
\end{Notation}

\begin{prop}\label{Fermatcaseextensionproperty}
In the Fermat case,
if $u$ is a locally convex function on $\partial \Delta_\lambda^\vee$, which is invariant under the permutation group. Then $u$ satisfies the extension property.
\end{prop}

\begin{proof}
We need to prove the characterisation in Prop. \ref{extensionpropertycharacterisation}.
Without loss of generality $L_\lambda(x)$ is achieved by $\langle m^0, x\rangle +\lambda(m^0)$. We need to find $p\in \Delta$, such that
$
u(y)- u(x) \geq \langle p, y-x\rangle.
$
For this we study the gradient of the function $u_{m^0}$ on the various $w$-charts.

First, notice for $x', y'$ on the face $\{ L_\lambda= \langle m^0, \rangle +\lambda(m^0)   \}$, namely the convex hull of $w_1, \ldots w_{n+1}$, the vector $y'-x'$ is parallel to the face, and by convexity of $u_{m^0}$ the directional derivative $\nabla u_{m^0}\cdot (y'-x')$ is monotone along the path from $x'$ to $y'$, so must be maximized at $y'$. In particular we consider such line segments on the face parallel to $w_i-w_j$ for $i, j\geq 1$. By the discrete symmetry, $\nabla u_{m^0}\cdot (w_j-w_i)$ must be zero on the plane of reflection bisecting the face. Thus for $i, j\geq 1$, $i\neq j$, the subset of the face
\[
\{ \nabla u_{m^0} \cdot w_i \geq   \nabla u_{m^0} \cdot w_j \} \cap \{ L_\lambda= \langle m^0, \rangle +\lambda(m^0)   \}
\]
agrees exactly with the half of the face containing $w_i$. Therefore the subset of face 
\[
\{ \nabla u_{m^0} \cdot w_i \geq   \nabla u_{m^0} \cdot w_j  , \forall j\geq 1  \}
\]
is exactly the intersection of  $\text{Star}(w_i)$ with the face. Without loss of generality $x$ lies in $\text{Star}(w_1)$.

We follow the notation in the proof of Prop. \ref{extensionimplieslocalconvexity}. In the $w_1$-chart, denote the gradient of $u_{m^0}$ as $\vec{p}$, so that for $y$ in the $w_1$-chart,
\[
u_{m^0} (y)- u_{m^0} (x)\geq \vec{p}\cdot (y-x)_{w_1}.
\]
A priori $\vec{p}$ lives in $M_\R/ \R m^0$. We lift $\vec{p}$ to $M_\R$ by demanding $\langle \vec{p}, w_1 \rangle=0$, so by the above discussion
$
\langle \vec{p},  w_i \rangle \leq 0 $ for $ i\geq 1. 
$
Define $p= \vec{p}+m_0$, then $\langle p, w_i\rangle \leq 1$ for all $i\geq 1$.
We regard $p\in M_\R$ as the gradient of $u$ at $x$, and write $p= \nabla u$ as a function of $x$. This construction can be made on other faces as well, and on the intersection of two faces the definitions are compatible.

We claim $p\in \Delta$: it suffices to show $\langle p, w_0\rangle \leq 1$. Notice $w_0= -\sum_1^{n+1} w_i=  \sum_{i=2}^{n+1} (w_1-w_i) -(n+1)  w_1 $. 
Consider the line segment in the face joining $x$ to the boundary of the face in the direction $ \sum_{i=2}^{n+1} (w_1-w_i)$, which stays inside $\text{Star}(w_1)$, and along which $\nabla u \cdot \sum_{i=2}^{n+1} (w_1-w_i)$ increases, or equivalently $\langle \nabla u, w_0\rangle$ increases. But the boundary of the face $\{ L_\lambda= \langle m^0, \rangle +\lambda(m^0)   \}$ lies also on a different face, and we can use the information from this new face to deduce $\langle \nabla u, w_0\rangle \leq 1$ there.

By construction for $y$ in the $w_1$-chart, 
\[
u(y)- u(x)\geq  \langle \vec{p}(x), (y-x)_{w_1}\rangle + \langle m_0, y-x\rangle = \langle \nabla u(x),  y-x  \rangle .
\]
We claim that in fact $u(y)-u(x)\geq \langle \nabla u(x), y-x\rangle $ holds for all $y\in \partial \Delta_\lambda^\vee$. We are left to check for $y$ on the face $\{ L_\lambda= \langle m^1, \rangle+ \lambda(m^1)  \}$, namely the complement of the $w_1$-chart. Consider the $w_i$-chart for $i>1$. We can write according to the decomposition $N_\R= (m^1)^\perp \oplus \R w_i$, that
\[
y-x= (y-x)_{w_i, m^1} + \langle m^1, y-x\rangle w_i.
\]
By local convexity, in the $w_i$-chart $u_{m^1}$ is convex, so there is some $\vec{p}'$, such that for any $y'$ in the $w_i$-chart
\[
u_{m^1}(y')- u_{m^1}(x) \geq \vec{p}' \cdot (y'-x)_{w_i, m^1}.
\] 
But a gradient vector of $u_{m^1}$ at $x$ is $\nabla u(x)- m^1$, so we may take $\vec{p}'= \nabla u(x)- m^1$. Thus
\[
u(y)-u(x)\geq \vec{p}' \cdot (y-x)_{w_i, m^1} + \langle m^1, y-x\rangle = \langle \nabla u(x), y-x\rangle - \langle \vec{p}' , w_i \rangle \langle m^1, y-x\rangle .
\]
Now $\langle m^1, y-x\rangle \geq 0$ as in the proof of Prop. \ref{extensionimplieslocalconvexity}, and $ \langle \vec{p}' , w_i \rangle  \leq 0$ by $\nabla u\in \Delta$. This implies $u(y)-u(x)\geq \langle \nabla u(x), y-x\rangle$ as required.

We have verified the characterisation in Prop. \ref{extensionpropertycharacterisation}, hence the extension property.	
\end{proof}

The proof above contains some additional information about the gradients.

\begin{cor}\label{extensiontranslationproperty}
In the region $\text{Star}(w_i)+ \R_{\geq 0} w_i \subset N_\R$, the directional derivative of the canonical extension $u= u^{**}$ satisfies
$
\langle w_i, \nabla u\rangle = 1.
$
In particular, in this region, for any $m$ with $\langle m, w_i \rangle =1$, the function $u_m=u-m$ is constant upon translation in the $w_i$-direction.

\end{cor}

\begin{proof}
By Remark \ref{pisgradient}, the $\nabla u$ introduced in the above proof is actually the gradient of the extension $u$ over $N_\R$. By the proof above, we know $\langle \nabla u, w_i \rangle=1 $ on $\text{Star}(w_i)\subset \partial \Delta_\lambda^\vee$. This directional derivative can only increase as $x\in N_\R$ moves in the $w_1$-direction. But $\nabla u \in \Delta$ on $N_\R$ since the extension is admissible, so  $\langle \nabla u, w_i \rangle\leq 1 $ everywhere, hence the claim.
\end{proof}

For later use, we define the notion of real MA equation in the Fermat case.

\begin{Def}\label{realMAFermat}
Let $u$ be a locally convex function  on $\partial \Delta_\lambda^\vee$ invariant under the discrete symmetry.  Then $u$ is called an Aleksandrov solution of the \emph{real MA equation} on $\partial \Delta_\lambda^\vee\setminus Sing$ if
\begin{itemize}
\item
On the interior of any top dimensional face of $\partial \Delta_\lambda^\vee$, in a set of standard local affine coordinates $x^{m_1}, \ldots x^{m_n}$ with $dx^{m_1}\wedge dx^{m_2}\ldots dx^{m_n}$ equal to the standard volume form $d\mu_\infty$, the function $u$ satisfies $MA(u)=d\mu_\infty$ in the Aleksandrov sense.

\item
On $\text{Star}(w)\subset U_w^\infty \cap \partial \Delta_\lambda^\vee$, we use the standard affine coordinates $x^{m_1}, \ldots x^{m_n}$ associated to the $U_w^\infty\cap \partial \Delta_\lambda^\vee$ chart. We demand for any vertex $m$ of $\Delta$ with $\langle m, w\rangle$, the local function $u_m=u-m$ satisfies $MA(u_m)=d\mu_\infty$ in the Aleksandrov sense.
\end{itemize}
Schematically we write $MA(u)= d\mu_\infty$.
\end{Def}

\begin{rmk}
Notice that the definition is compatible on overlapping charts because the transition functions lie in $SL(n,\Z) \ltimes \R^n $. On the locus $Sing\subset \partial \Delta_\lambda^\vee$ we make no definition.	
\end{rmk}

\section{Estimates on the K\"ahler potential}\label{Estimatesonpotential}

This section is concerned with estimating the K\"ahler potential on the degenerating hypersurfaces $X_s$ in the Fermat family. The expectation that the potentials converge in the $s\to + \infty$ limit to a solution of a real MA equation, motivates us to produce local convex functions by taking average of local K\"ahler potentials.  Convex functions have better a priori regularity than psh functions: a Lipschitz bound is automatic.  These arguments work for general K\"ahler potentials, without using the complex MA equation. The main difficulty is then to show that for the Calabi-Yau metric, the local potentials are $C^0$-close to their averaging convex functions at least in the generic region; equivalently the local potentials have small local oscillations. This part relies on the method of Kolodziej as outlined in section \ref{Toolsfrompsh}, and a key ingredient is an improved uniform Skoda inequality.

Most arguments apply to more general contexts, and the only reason we restrict to the Fermat family of hypersurfaces is to use the extension property, which enables us to patch up the local convex functions into a global regularisation of the original K\"ahler potential.

\subsection{Harnack inequality}

Consider a general possibly singular K\"ahler potential $\varphi\in PSH(X_s, s^{-1}\omega_{FS})$ on $X_s$, normalised to $\sup_{X_s} \varphi=0$.
We think of $\varphi$ equivalently as a collection of local potentials $\{\varphi_0, \varphi_m\}$ as in section \ref{Kahlerianpolarisation}. In the region $U_w^s\subset X_s$,  we can find $m\in \Delta_\Z$ with $\langle m, w\rangle=1$ and $\C^*$-coordinates $z^{m_1}, \ldots z^{m_n}$ as  in section \ref{Complexstructure}. Recall $d\mu_s$ is the normalised canonical measure induced by the holomorphic volume form.

\begin{Notation}
Denote $X^{toric}_s$ as the union of all the toric regions $ U_{w,\delta}^s$ for various choices of $m$ and $w$. It is tacitly understood that slightly shrinked domains correspond to a slightly larger choice of $\delta$, and we shall abusively use the same notation for shrinked domains.	
\end{Notation}

\begin{prop}\label{Harnacktypeinequality}
	(Harnack type inequality) Suppose $\varphi\in PSH(X_s, s^{-1}\omega_{FS})$ with $\sup_{X_s} \varphi=0$. Then the average integral
	\[
	\dashint_{X^{toric}_s } |\varphi| d\mu_s \leq C.
	\]
	
\end{prop}

\begin{proof}
	(\cf proof of Prop. 3.1 in \cite{Blockilecture})
	Consider the local potentials $\phi=\varphi_m$ on various coordinate charts in section \ref{Complexstructure}, both of the toric type and of the boundary type. The charts can be chosen so that the Lebesgue measures thereof are uniformly equivalent to $d\mu_s$ up to a scaling factor. We have $|\phi- \varphi|\leq C$ uniformly on charts. Suppose a coordinate ball $B(p, 3R)$ is contained in (the universal cover of) the local chart. Since $\phi$ is psh and $\phi-C\leq 0$, for $z\in B(p,R)$,
	\[
	\phi(y)-C \leq \dashint_{B(y,2R)} (\phi-C) \lesssim   \dashint_{B(p,R)} (\phi-C),
	\]
	hence
	\[
	\dashint_{B(p,R)} |\varphi |\lesssim 1+ \inf_{B(p, R)} (- \varphi). 
	\]
	To deduce the global version of the Harnack type inequality we need a transitivity property, namely we can connect the chart containing the maximum point of $\varphi$ to any of the toric charts in $X_s^{toric}$ via a chain of $O(1)$ number of charts, such that  $\inf_{B(p,R)} |\varphi|$ on charts increase by only $O(1)$ in each step. This last fact is because we can choose the chains of successive charts $B(p_i, 5R_i)$ such that the measure of the overlap occupies a nontrivial portion of the previous chart:
	\[
	| B(p_i,R_i)\cap B(p_{i+1}, R_{i+1}) | \gtrsim  \frac{1}{10} | B(p_i,R_i) |,
	\]
	which would force \[
	\begin{split}
	\inf _{ B(p_{i+1},R_{i+1}) } |\varphi| \leq & \inf _{ B(p_{i+1},R_{i+1})\cap  B(p_i,R_i) } |\varphi| \leq \dashint_{ B(p_{i+1},R_{i+1})\cap  B(p_i,R_i) } |\varphi| 
	\\
	\lesssim  & \dashint_{  B(p_i,R_i) } |\varphi|  \lesssim 1+ \inf _{ B(p_{i},R_{i}) } |\varphi|.
	\end{split}
	\]
	
\end{proof}

\begin{rmk}
Notice this transitivity argument allows us to move from boundary type charts into toric charts, but not conversely, because the measure is much larger on toric charts.
\end{rmk}

\subsection{Local potentials: convexity }\label{Localpotentialconvexity}

We continue with a general $\varphi\in PSH(X_s, s^{-1}\omega_{FS})$ normalised to $\sup_{X_s} \varphi=0$, whose local potentials are $\{\varphi_0, \varphi_m\}$.
A simple obeservation is:

\begin{lem}\label{pshconvexity}
Let $\Phi$ be any psh function on the open subset of $\{   1<|\zeta_i|<\Lambda, i=1,\ldots n \}\subset (\C^*)^n$. Then the $T^n$-invariant function \[
\bar{\Phi}(\log |\zeta_1|, \ldots, \log |\zeta_n|)= \frac{1}{(2\pi)^n} \int_{T^n} \Phi( |\zeta_1|e^{i\theta_1}, \ldots |\zeta_n|e^{i\theta_n}) d\theta_1\ldots d\theta_n
\]
is a convex function in the variables $x_1=\log |\zeta_1|, \ldots, x_n=\log |\zeta_n|$. 
\end{lem}

\begin{proof}
Since the $T^n$-action on $(\C^*)^n$ is holomorphic, $\Phi(\zeta_1 e^{i\theta_1}, \ldots \zeta_n e^{i\theta_n})$ is psh in $\zeta$ for any choice of $\theta_i$, so the average function $\bar{\Phi}$ is also psh. Any $T^n$-invariant psh function must be convex in the log coordinates, because of the formula
\[
\sqrt{-1}\partial \bar{\partial} \bar{\Phi}= \frac{1}{4} \sum \frac{\partial^2 \bar{\Phi}}{  \partial x_i \partial x_j  } \sqrt{-1} d\log \zeta_i \wedge d\overline{\log \zeta_j} \geq 0.
\]
\end{proof}


In the region $U_w^s\subset X_s$,  we can find $m\in \Delta_\Z$ with $\langle m, w\rangle=1$ and $\C^*$-coordinates $z^{m_1}, \ldots z^{m_n}$ as  in section \ref{Complexstructure}, and consider the local potential $\phi=\varphi_m$. Denote $x^{m_i}=\frac{\log |z^{m_i}|}{s}$. We produce the local average function
\begin{equation}
\bar{\phi}(x^{m_1},\ldots x^{m_n} )= \frac{1}{(2\pi)^n} \int_{T^n} \phi( |z^{m_1}|e^{i\theta_1}, \ldots |z^{m_n}|e^{i\theta_n}) d\theta_1\ldots d\theta_n.
\end{equation}

\begin{prop}\label{Lipschitzboundlocalpotential}
In the chart $ U^s_w $ the average function $\bar{\phi}$ is convex, 
 and on the shrinked chart $ U^s_{w,\delta}  $ it
has a Lipschitz bound:
\begin{equation}\label{Lipschitzbound}
|\bar{\phi}|\leq C, \quad 
|\bar{\phi}(x)- \bar{\phi}(x')| \leq C|x-x'|.
\end{equation}
\end{prop}

\begin{proof}
By Lemma \ref{pshconvexity},   $\bar{\phi}$ is convex, and by Prop. \ref{Harnacktypeinequality} it has an $L^1$ bound in the $x^{m_i}$ coordinates:
\[
\int |\bar{\phi}| dx^{m_1}\ldots dx^{m_n} \leq C.
\]
Clearly $\bar{\phi}$ is also bounded above, so for the argument we may pretend $\bar{\phi}\leq 0$ upon shifting by a bounded constant.

We claim $\bar{\phi}(x)$ is bounded from below for $x$ in a shrinked interior region. The ball $B(x,2r)$ is contained in the coordinate chart, with $r$ bounded below by a positive constant. For $y$ in the annulus $B(x,2r)\setminus B(x, r)$, we have $2\bar{\phi}( \frac{x+y}{2} )\leq \bar{\phi}(x ) + \bar{\phi}(y )$,	so upon integration
\[
\int | \bar{\phi}| \gtrsim \int 2|\bar{\phi}( \frac{x+y}{2} )| dy \geq \int |\bar{\phi}(x )| + |\bar{\phi}(y )|dy,
\]	
which bounds $|\bar{\phi}(x)|$. Thus on a slightly shrinked $x$-domain  the oscillation is bounded:
\[
\text{osc } \bar{\phi}=(\sup-\inf )\bar{\phi} \leq C,
\]
and the Lipschitz bound follows again by convexity.
\end{proof}

\begin{rmk}
	We discuss some intuition about \emph{log scales}. Let $P\in X_s$ lie in $ U_{w,\delta}^s $, then a log scale $|z^{m_i}|\sim |z^{m_i}(P)|$ around $P$ refers to the subregion
	\[
	\{  \frac{1}{2} |z^{m_i}(P)|  \lesssim|z^{m_i}| \lesssim 2|z^{m_i}(P)|  , \quad 1\leq i\leq n    \}.
	\] 
	Now $\log |z^{m_i}|$ vary by order $O(s)$ within $ U_{w,\delta}^s $, so there are an enormous number of log scales. The long range behaviour of $X_s$ is similar to $(\C^*)^n$, with half of the dimensions compactified into $T^n$. On the other hand, over one log scale $X_s$ behaves qualitatively like the unit disc in $\C^n$. The concept of \emph{local oscillation} of a function refers to the oscillation within one log scale. In particular the Lipschitz bound (\ref{Lipschitzbound}) implies a local oscillation bound
	\[
	\text{osc}_{ |z^{m_i}|\sim |z^{m_i}(P)| } \bar{\phi} \leq Cs^{-1}.
	\]
\end{rmk}

\subsection{Local potentials: plurisubharmonicity}

The following lemma is a special case of the principle that for a subharmonic function, the standard mean value inequality has interesting strengthenings if there is more information about microscopic averages.

\begin{lem}
Let $\Phi$ be a subharmonic function on $B_2^n\times T^k= B_2 \times \R^k/ \epsilon \Z^k$ equipped with the Euclidean metric $g=\sum_1^n dx_i^2+ \sum_1^k dy_j^2$, where $0< \epsilon \ll 1$. Let $v$ be the averaging function of $\Phi$ over the $T^k$ fibres. Assume $\dashint |\Phi| \lesssim 1$ and a Lipschitz bound $\text{Lip}(v) \lesssim 1$, then on $B_1\times T^k$ we have
$\Phi\leq v+ C\epsilon^{1/2}$. 
\end{lem}

\begin{proof}
(courtesy of W. Feldman) By passing to the universal cover $B_2\times \R^k$, the standard mean value inequality implies
\[
\sup_{ B_{3/2}\times T^k   } \Phi \lesssim \dashint |\Phi| \lesssim 1.
\]
Let $p\in B_1\times T^k$, which lifts to a point $p$ in $B_1\times \R^k$. Consider the Euclidean ball $B_g(p, \epsilon R)\subset B_{3/2}\times \R^k$, where $R\gg 1$ is a parameter to be chosen. Then by the mean value inequality,
\[
\Phi(p) \leq \dashint_{ B_g(p, \epsilon R)  } \Phi .
\] 
Define the subset $E\subset B_g(\epsilon R)$ as the union of all interior lattice cubes, then
\[
B_g(p, \epsilon R  )\setminus E \subset B_g ( p, \epsilon R  )\setminus B_g(p, \epsilon (R-C)) ,
\]
and by the lattice periodicity of $\Phi$ we have  $\int_E \Phi= \int_E v$. By partitioning 
the integral $\int_{ B_g ( p,\epsilon  ) }\Phi$ into the contributions from $E$ and $B_g(p, \epsilon R  )\setminus E $,
\[
\dashint_{ B_g(p, \epsilon R)  } \Phi 
\leq  \dashint_{ B_g(p, \epsilon R)     } v +  CR^{-1} \sup_{B_g(p,\epsilon R)  } (\Phi-v) \leq \dashint_{ B_g(p, \epsilon R)     } v +  CR^{-1} .
\]
By the Lipschitz bound of $v$, the RHS is bounded above by
\[
v(p) +  \text{Lip}(v) \epsilon R+  CR^{-1} \leq v(p)+ C(\epsilon R+ R^{-1} ).
\]
Choosing $R= \epsilon^{-1/2}$ gives $\Phi(p)\leq v(p)+ C\epsilon^{1/2}$.
\end{proof}

Back to the setting of Prop. \ref{Lipschitzboundlocalpotential},

\begin{cor}\label{potentialupperboundedbyconvex}
(Local potential upper bound)
On $ U^s_{w,\delta}  $, then
$
\phi- \bar{\phi} \leq Cs^{-1/2}.
$
\end{cor}

\begin{proof}
The psh property of $\phi$ implies subharmonicity.
By the Harnack inequality in Prop. \ref{Harnacktypeinequality} the average $L^1$-integral is bounded, and by Prop. \ref{Lipschitzboundlocalpotential} there is a Lipschitz bound on the local average function $\bar{\phi}$.
\end{proof}

\begin{cor}\label{localL1oscillationI}
(Local $L^1$-oscillation bound) Over one log scale inside $U_{w,\delta}^s$,
\[
\dashint_{ |z^{m_i}|\sim |z^{m_i}(P)|  } |\phi- \bar{\phi}| d\mu_s \leq C s^{-1/2}.
\]
\end{cor}

\begin{proof}
Recall the local oscillation of $\bar{\phi}$ in one log scale is $O(s^{-1})$. Since the local sup of $\phi$ differs from the local average of $\phi$ by $O(s^{-1/2})$, the local $L^1$-oscillation is likewise bounded by $O(s^{-1/2})$.
\end{proof}

\begin{rmk}
	The $s$-dependence is probably not optimal.
\end{rmk}



We now seek a local $L^1$-oscillation bound on the charts of boundary type $U_P$ (\cf Remark \ref{chartsofboundarytype}). The idea is that any chart of boundary type overlaps with some chart of toric type
in an annulus region, where the $L^1$-oscillation bound is already known. It would be enough to transfer the $L^1$-oscillation bound from the annulus to the deep interior of the chart.

\begin{lem}
Let $\Phi$ be a psh function on the  $\{  |z_i|\leq 4, \forall i \}\subset \C^n$. Then
\[
\dashint_{ B_1 } |\Phi| \lesssim  \dashint_{ \{ 1<|z_i|<4,\forall i   \}  } |\Phi|   .
\]
\end{lem}

\begin{proof}
We induct on dimension. For $n=1$, the unit ball is already enclosed by an annulus, so $\sup_{B(1)} \Phi$ is bounded above, and the mean value property applied to all balls $B(p,2)$ with $1<|p|\leq 2$ gives a lower bound on $\dashint_{B(1) } \Phi$. Thus the $L^1$-bound in $B(1)$ is clear.

For general $n$, notice by induction we can bound for each $i\leq n$,
\[
\dashint_{ \{ 1< |z_i|<4 , |z_j|< 4, \forall j\neq i  \}  } |\Phi|  \lesssim  \dashint_{ \{ 1<|z_j|<4,\forall j   \}  } |\Phi|,
\]
so $\Phi$ is controlled in $L^1$ on an annulus enclosing $B(1)$, and we can bound $\dashint_{B(1)} |\Phi|$ similar to the $n=1$ case.
\end{proof}

\begin{cor}\label{localL1oscillationII}
(Local $L^1$-oscillation bound II) In the chart of boundary type $U_P$, the local potential $\phi$ satisfies
\[
\dashint_{U_P} |\phi- \dashint_{U_P} \phi| d\mu_s \leq C s^{-1/2}.
\]
\end{cor}

\subsection{Locally convex function}\label{Locallyconvexfunctionsection}

In section \ref{Localpotentialconvexity}  we produced a collection of local average functions $\bar{\phi}=\bar{\phi}_{m,w}$ on $U_w^s$ corresponding to various choices of $w$ and $m$ with $\langle m,w\rangle=1$. But the local coordinates $x^{m_1}, \ldots x^{m_n}$ are naturally interpreted also as coordinates on $\partial \Delta_\lambda^\vee$ (\cf section \ref{Piecewiseaffinestructure}), so $\bar{\phi}_{m,w}$ can be alternatively viewed as a collection of convex functions on the charts $U_w^\infty\cap \partial \Delta_\lambda^\vee$ of $\partial \Delta_\lambda^\vee$. (Notice these local functions are defined without the need to shrink the domain to $U^\infty_{w,\delta}$).

The intuition is that up to $C^0$-small error, the differences of these local functions agree with the cocycle $\{ m-m' \}$, or equivalently, up to some $C^0$-small fuzziness $\bar{\phi}_{ m, w}+ \langle m, x \rangle$ 
glue to a \emph{locally convex function} on $\partial \Delta_\lambda^\vee$ in the sense of Definition \ref{locallyconvexdef}. The more precise statement is

\begin{lem}
On overlapping charts of $\partial \Delta_\lambda^\vee$, 
\[
| \bar{\phi}_{ m, w}  - \bar{\phi}_{ m', w'} +(m-m') |
\leq C s^{-1/2}.
\]	
\end{lem}

\begin{proof}
Since we know the local $L^1$-oscillation estimate holds in every local region, in a log scale in $U^s_w$, not necessarily in the shrinked region $U^s_{w,\delta}$, 
\[
\dashint_{ |z^{m_i}|\sim |z^{m_i}(P)|  } |\varphi_m - \dashint \varphi_m | d\mu_s \leq Cs^{-1/2}.
\]
Since $\bar{\phi}_{m,w}$ is convex, a local $L^1$-bound implies a local $L^\infty$-bound in a slightly shrinked region, so in the log scale,
\[
| \bar{\phi}_{m,w}- \dashint \varphi_m |\leq Cs^{-1/2}. 
\]
Likewise for $\bar{\phi}_{m', w'}$.  By definition the local potentials differ by 
\[
\varphi_m- \varphi_{m'}= \langle m'-m, \text{Log}_s(z) \rangle
\]
Notice that for a given point $P$ on $\partial \Delta_\lambda^\vee$, the log scales on $U^s_w$ and $U^s_{w'}$ around $P$ have a nontrivial percentage of overlapping measure. Thus
\[
\begin{split}
|\bar{ \phi}_{m,w}- \bar{\phi}_{m', w'} +(m-m') |& \lesssim s^{-1/2} + |  \dashint \varphi_m - \dashint \varphi_{m'} +(m-m') |
\\
& \lesssim s^{-1/2} + \dashint_{overlap} |\varphi_m- \varphi_{m'} + (m-m') |
\\
& \lesssim s^{-1/2}.
\end{split}
\]
\end{proof}

\begin{rmk}
The tropical version $U^\infty_w$ of $U^s_w$ is in general larger than $U^\infty_w\cap \partial \Delta_\lambda^\vee$; it typically contains also some subset stretching to infinity along the $w$-direction. If we regard   $\bar{\phi}_{m,w}$ as local functions on $\mathcal{A}_\lambda^\infty$ instead of $\partial \Delta_\lambda^\vee$, then there is a delicate issue.
 The Lemma above \emph{does not} imply that $\bar{\phi}_{m,w}+ \langle m,x\rangle$ for various choices of $m,w$ glue approximately on overlapping regions far from $\partial \Delta_\lambda^\vee$. The problem is that such overlapping regions have too small measure, which breaks down the proof. 
\end{rmk}

\subsection{Legendre transform, extension,  regularisation}\label{Legendretransform}

We restrict to the Fermat case, and consider a general $\varphi\in PSH(X_s, s^{-1}\omega_{FS})$ with $\sup_{X_s} \varphi=0$, invariant under the symmetric group permuting the monomials $Z_0^{n+2}, \ldots, Z_{n+1}^{n+2}$.
The goal of this section is to canonically patch together the local convex functions in section \ref{Locallyconvexfunctionsection} approximately to produce a convex admissible function on $N_\R= \R^{n+1}$.
We will then induce a potential $\psi\in PSH(X_s, s^{-1}\omega_{FS} )\cap C^0$ which is a \emph{regularisation} of $\varphi$ in the sense that it enjoys better a priori bounds than $\varphi$.

\begin{prop}\label{regularisationFermat}
There is an admissible convex function $u$ on $N_\R$, such that on $U_w^\infty\cap \partial \Delta_\lambda^\vee$,
\begin{equation}
|u - ( \bar{\phi}_{m,w} + m      )|\leq Cs^{-1/2}.
\end{equation}
\end{prop}

\begin{proof}
The idea is to regard $ \bar{\phi}_{ m, w}+ \langle m, x \rangle$ as approximately defining a locally convex function on $\partial \Delta_\lambda^\vee$ in the sense of Def. \ref{locallyconvexdef}, and then the problem is essentially to prove an effective version of the extension property (\cf Prop. \ref{Fermatcaseextensionproperty}). We will outline the main modifications.

We will produce $u$ by mimicking the Legendre duality construction in Prop. \ref{extensionpropertycharacterisation}.
For $p\in \Delta$, define 
\[
u^*(p)= \sup_{ x\in \partial \Delta_\lambda^\vee } \{    \langle x,p\rangle - ( \bar{ \phi}_{ m, w}+ \langle m, x \rangle    )       \},
\]
where it is tacitly understood that $ \bar{\phi}_{ m, w}+ \langle m, x \rangle  $ is defined only over $\partial \Delta_\lambda^\vee\cap U^\infty_w$, and the sup is taken over all choices of $m, w$ whenever $\bar{\phi}_{m,w}$ is defined. Since $\bar{\phi}_{m,w}$ are uniformly bounded on $\partial \Delta_\lambda^\vee$, we see $\norm{u^*}_{C^0(\Delta)} \leq C$. We then define a convex function $u$ on $N_\R$ by another Legendre transform
\[
u (x)= \sup_{p\in \Delta}  \{ \langle p, x\rangle - u^*(p)   \},
\]
which is admissible because $u^*$ is bounded. By the same reasoning in Prop. \ref{extensionpropertycharacterisation}, on $\partial \Delta_\lambda^\vee\cap U^\infty_w$, 
\[
u(x) \leq \bar{\phi}_{m, w}+ \langle m, x \rangle + Cs^{-1/2}  .
\]
We are only left to show
\[
u(x)\geq \bar{\phi}_{m, w}+ \langle m, x \rangle - Cs^{-1/2} ,
\]
which amounts to showing that there exists $p\in \Delta$, such that for any $y\in \partial \Delta_\lambda^\vee$, 
\[
\bar{\phi}_{m', w'}(y)+ \langle m', y \rangle \geq \bar{\phi}_{m, w}(x)+ \langle m, x \rangle +\langle p, y-x\rangle- Cs^{-1/2}.
\]

Notice our setting enjoys the discrete symmetry.
This last step is the effective version of Prop. \ref{Fermatcaseextensionproperty}, and the proof is basically the same. 
\end{proof}

By construction $u$ has a number of additional properties:
\begin{cor}\label{canonicalextensiontranslation}
The canonical extension $u$ satisfies an a priori Lipschitz bound
\begin{equation}\label{Lipschitzboundonu}
\begin{cases}
|u- \max_m \langle m, x\rangle| \leq C, \quad \forall x\in N_\R,\\
|u (x) - u(x') | \leq C|x-x'|, \quad \forall x, x'\in N_\R.
\end{cases}
\end{equation}
Morever, in the region $\text{Star}(w)+\R_{\geq 0} w\subset N_\R$, for any $m$ with $\langle m, w\rangle=1$, the function $u_{m}=u-m$ is constant upon translation in the $w$-direction.

\end{cor}

\begin{proof}
The first inequality is because the Legendre transform $u^*(p)$ is bounded on $\Delta$ as in the above proof, and the second is because $\nabla u\in \Delta$. The morever statement is essentially identical to Cor. \ref{extensiontranslationproperty}.
\end{proof}

By a small variant of Prop. \ref{admissibleimpliesextension}, when we pullback the admissible convex functions $u$ via $\text{Log}_s$, we obtain a torus invariant K\"ahler current on $(\mathbb{P}_\Delta, s^{-1}[\Delta])$ with continuous local potentials. In details, we write $\psi_0= u \circ \text{Log}_s$, and define
\begin{equation}
\begin{cases}
\psi= \psi_0 -  \frac{ (n+2)}{2s} \log (\sum_{m'\in vertex(\Delta) } e^{ \frac{2}{n+2} \langle m', \text{Log}(z) \rangle }  )  ,
\\
\psi_{m}= \psi_0 - \langle m, \text{Log}_s(z) \rangle = u_m \circ \text{Log}_s .
\end{cases}
\end{equation}
By construction $\psi\in PSH(\mathbb{P}_\Delta, s^{-1} \omega_{FS}) \cap C^0$, and $\psi_0, \psi_m$ are the local potentials of $\psi$ (\cf  (\ref{localpotentialscomplexsetting})). By Cor. \ref{canonicalextensiontranslation}, $\norm{\psi}_{C^0} \leq C$, and $\psi$ inherits the Lipschitz bound from $u$.
By a slight abuse of notation, the restriction to $X_s$ will still be denoted as $\psi\in PSH(X_s,s^{-1}\omega_{FS} )\cap C^0$.
We think of $\psi$ as a regularisation of $\varphi$.

\begin{rmk}
As explained in section \ref{Algebraicmetrics}, on toric manifolds the Legendre transform arises from a limiting version of approximation by algebraic metrics, which in turn is a more standard way to regularise an arbitrary K\"ahler potential. Now $X_s$ is not a toric manifold, but the toric symmetry holds approximately in generic regions, which motivates us to take the Legendre transform as a replacement of algebraic regularisation. 
\end{rmk}

We now specify some subregions on $X_s^{toric}$ with coordinate descriptions. These are intimately related to $\partial \Delta_\lambda^\vee\setminus Sing$, which is covered by the stars of the vertices and the interior of the top dimensional faces (\cf section \ref{Fermatcaseextensionpropertysection}).

\begin{Notation}(Star type regions on $X_s$)
On the region $U^s_w\subset X_s$, recall the coodinates $z^{m_i}$ and regard $x^{m_i}= s^{-1}\log |z^{m_i}|$ as local coordinates also on $U^\infty_w\cap \partial \Delta_\lambda^\vee$. Let 
 $U^{s,*}_w\subset U^s_{w,\delta}$ be the subset where the $x^{m_i}$ coordinates correspond to points in $\text{Star}(w)\subset \partial \Delta_\lambda^\vee$. The tropical analogue of  $U^{s,*}_w$ is 
  $(\text{Star}(w)+\R_{\geq 0} w )\cap A_\lambda^\infty$.
\end{Notation}

\begin{Notation}
(Face type regions on $X_s$) Consider a slightly shrinked subset of the interior of a given top dimensional face of $\partial \Delta_\lambda^\vee$. This can be regarded as a subset of $U_{w,\delta}^\infty\cap \partial \Delta_\lambda^\vee$, where we regard $x^{m_i}= s^{-1}\log |z^{m_i}|$ as local affine coordinates. Let $U_w^{s, face}\subset U_w^s$ be the subset where the $x^{m_i}$ coordinates correspond to points in this shrinked face. The tropical analogue of $U_w^{s, face}\subset U_w^s$ is the shrinked face.

\end{Notation}

The intuition is that when $z\in X_s$ has $\text{Log}_s$ image close to $\partial \Delta_\lambda^\vee$, or if this image approaches infinity in specific directions, then $\varphi-\psi$ is bounded above by a very small number:

\begin{prop}\label{LocalpotentialupperboundII}
(Local potential upper bound)
\begin{itemize}
\item
Inside $U^{s,*}_w\subset X_s$, for $\langle m,w\rangle=1$, the local potentials satisfy
$
\varphi_m - \psi_m \leq  Cs^{-1/2},
$
or equivalently $\varphi-\psi \leq Cs^{-1/2}$.

\item
Inside $U^{s,face}_w$, the local potentials satisfies $\varphi_0- \psi_0 \leq Cs^{-1/2}$, or equivalently $\varphi-\psi\leq Cs^{-1/2}$.

\end{itemize}

\end{prop}

\begin{proof}
In the star type region case,
by Cor. \ref{potentialupperboundedbyconvex}, we have the upper bound $\varphi_m- \bar{\phi}_{m,w} \leq Cs^{-1/2}$. 
By Prop. \ref{regularisationFermat} and Cor. \ref{canonicalextensiontranslation}, in $U^{s,*}_w$ we can replace $\bar{\phi}_{m,w}$ by $\psi_m$ up to an error bounded by $Cs^{-1/2}$, hence the claim. The face type region follows the same argument, without the translational invariance statement of Cor. \ref{canonicalextensiontranslation}. 
\end{proof}

\subsection{Improved Skoda inequality}

Recall the local $L^1$-oscillation bounds in both toric and boundary type regions, from Cor. \ref{localL1oscillationI} and \ref{localL1oscillationII}. Consequently,

\begin{lem}
	(Local Skoda estimate) Consider any $\varphi \in PSH(X_s , s^{-1} \omega_{FS})$ normalised to $\sup_{X_s} \varphi=0$.
	There are uniform positive constants $\alpha$, $C$, such that the local potentials $\phi$ satisfy
\begin{itemize}
\item In a log scale in the toric region,
\[
\dashint_{ |z^{m_i}|\sim |z^{m_i}(P)|  } e^{- \alpha \sqrt{s} ( \phi- \dashint \phi   ) } d\mu_s \leq C.
\]
\item In a boundary type chart,
\[
\dashint_{ U_P } e^{- \alpha \sqrt{s} ( \phi- \dashint \phi   ) } d\mu_s \leq C.
\]
\end{itemize}

\end{lem}

\begin{proof}
Apply the standard Skoda inequality  (\cf Thm \ref{Skodabasicversion}) to the rescaled function $s^{1/2}(\phi- \dashint \phi )$.
\end{proof}

\begin{rmk}
The local average $\dashint \phi$ can be replaced by the local supremum using the mean value inequality.
\end{rmk}

\begin{cor}\label{globalSkoda}
	(global Skoda estimate) Consider any $\varphi \in PSH(X_s , s^{-1} \omega_{FS})$ normalised to $\sup_{X_s} \varphi=0$. There are uniform positive constants $\alpha$, $C$, such that
	\begin{equation}
	\dashint_{X_s}  e^{- \alpha  \varphi}  d\mu_s \leq C.
	\end{equation}
\end{cor}

\begin{proof}
By the local Skoda estimate and the Remark above, for both a log scale in the toric region, and a boundary type chart, the local average
\begin{equation}\label{LocalSkodapreperation}
\dashint e^{ \alpha \sqrt{s} (-\varphi+ \sup_{loc} \varphi ) } d\mu_s \leq C,	
\end{equation}
so in particular
$
\dashint e^{ \alpha  (-\varphi+ \sup_{loc} \varphi ) } \leq C. 
$
But we have already  achieved a $C^0$-bound on local average functions, and in particular a lower bound on local suprema. Thus
\[
\dashint e^{ -\alpha \varphi } d\mu_s \leq C e^{-\alpha  \sup_{loc} \varphi  } \leq C,
\]
or equivalently $\int_{loc} e^{ -\alpha \varphi } d\mu_s \leq C \int_{loc} d\mu_s$ for local integrals. To pass from this to the global Skoda estimate, we need to take a large collection of log scales and boundary type charts and sum over the estimates:
\[
\int_{X_s} e^{-\alpha \varphi} d\mu_s \leq C\sum \int_{loc} d\mu_s.
\]
The only problem is to ensure that the local charts can be chosen without substantially overcounting the measure. For points on $X_s$ whose $\text{Log}_s$ image is at $O(s^{-1})$ Euclidean distance to $\partial \Delta_\lambda^\vee$, it is easy to choose the charts so that each point is contained in $O(1)$ number of charts. Away from $\partial \Delta_\lambda^\vee$, the points deep inside the boundary type charts in general do not have this local finiteness property, but this is compensated by the fact that the measure $d\mu_s$ decays exponentially away from $\partial \Delta_\lambda^\vee$ (\cf (\ref{exponentialmeasuredecayawayfromtoric})). The conclusion is that
\[
 \sum \int_{loc} d\mu_s \leq C\int_{X_s} d\mu_s,
\]
whence the global Skoda estimate.
\end{proof}

We now specialize to the Fermat  case, and consider $\varphi\in PSH(X,s^{-1} \omega_{FS})$ normalised to $\sup_{X_s} \varphi=0$ with discrete symmetry, as in section \ref{Legendretransform}. The regularisation of $\varphi$ produced via Legendre transform is denoted as $\psi$.

\begin{thm}\label{improvedSkoda}
(Improved Skoda estimate) 
In the Fermat case above,  there are uniform constants $\alpha$, $C$, such that
\begin{equation}
\dashint_{X_s}  e^{- \alpha \sqrt{s} ( \varphi- \psi   ) } d\mu_s \leq C.
\end{equation}
\end{thm}

\begin{proof}
On either a log scale in the toric region, or a boundary type chart, we have by the local $L^1$-oscillation estimate and the mean value inequality that
\[
|\sup_{loc} \varphi_m - \dashint _{loc}\varphi_m|\leq Cs^{-1/2}, \quad |\sup_{loc} \psi_m- \dashint_{loc} \psi_m | \leq Cs^{-1/2}.
\]
Notice also the local averages of $\varphi_m$ and $\psi_m$ differ by $O(s^{-1/2})$, so 
\[
|\sup_{loc} \varphi_m- \sup_{loc} \psi_m |\leq Cs^{-1/2}, \quad |\sup_{loc} \varphi- \sup_{loc} \psi |\leq Cs^{-1/2}.
\]
Combined with (\ref{LocalSkodapreperation}),
\[
\dashint_{loc}  e^{- \alpha \sqrt{s} ( \varphi- \psi   ) } d\mu_s \leq C.
\] 
The summation argument as in the global Skoda estimate proves the claim.
\end{proof}

\begin{rmk}
This means $\phi-\psi$ can only fail to be bounded below by $Cs^{-1/2}$ on a set with exponentially small  probability measure. Notice we have not yet used the complex MA equation.
\end{rmk}

\subsection{$L^\infty$ and stability estimates for CY potentials}\label{Localoscillationestimate}

We finally impose the Calabi-Yau condition, and consider the CY potential $\varphi= \varphi_{CY,s}$ normalised to $\sup_{X_s} \varphi=0$, solving (\ref{CalabiYaucondition}):
\[
\omega_{CY,s}^n= ( s^{-1} \omega_{FS} + \sqrt{-1} \partial \bar{\partial} \varphi )^n= a_s s^{-n} d\mu_s.
\]

\begin{thm}
($L^\infty$-estimate) The Calabi-Yau potential $\varphi_{CY,s}$ satisfies the uniform $L^\infty$-estimate $\norm{\varphi_{CY,s}}_{L^\infty} \leq C$.
\end{thm}

\begin{proof}
We apply Kolodziej's estimate in Thm \ref{pluripotentialthm1}. The Skoda type inequality (\ref{Skodaassumption}) is verified in Cor. \ref{globalSkoda}, hence the $L^\infty$ estimate. 
\end{proof}

We now specialize to the Fermat case.
Clearly $\varphi$ is invariant under the discrete symmetry of the  hypersurface. Recall the regularisation is denoted as $\psi=\psi_{CY,s}$, coming from the double Legendre transform construction $u=u_{CY,s}$ (\cf section \ref{Legendretransform}).
The local potentials of $\varphi_{CY,s}$ and $\psi_{CY,s}$ are denoted $\varphi_m=\varphi_{CY,s,m}$ and $\psi_m= \psi_{CY,s,m}$ according to the same convention as (\ref{localpotentialscomplexsetting}).

\begin{thm}
In the Fermat case, 
there is a uniform stability estimate \begin{equation}
\varphi_{CY,s} - \psi_{CY,s} \geq -Cs^{-1/2} \log s  .   
\end{equation}
\end{thm}

\begin{proof}
We apply Cor. \ref{StabilityestimateKolodziej}. The Skoda estimate is verified in Cor. \ref{globalSkoda}. The improved Skoda estimate Thm. \ref{improvedSkoda} implies an exponential volume decay:
\[
\frac{ \int_{ \varphi- \psi \leq -t} \omega_\phi^n}{ \text{Vol}(X_s)  } \leq C e^{-\alpha t \sqrt{s} }, 
\]
hence there exists $c\gg 1$, such that for $t_0 = c s^{-1/2} \log s$,
\[
\left(\frac{ \int_{ \varphi- \psi \leq -t_0} \omega_\phi^n}{  \text{Vol}(X_s)  } \right)^{1/2n}  \leq C e^{-\alpha t_0 \sqrt{s}/2n } = C e^{ - \alpha c \log s/ 2n  } \leq C s^{-1/2}.
\]
Thm \ref{pluripotentialthm1} then implies $\varphi - \psi \geq -Cs^{-1/2} \log s  $ as required.
\end{proof}

\begin{rmk}
In the theorems above only an upper bound on the volume measure is actually needed. The intuition is that the Skoda inequality is already so close to an $L^\infty$ estimate, that a very tiny amount of extra assumptions are needed to conclude $L^\infty$-estimate.
\end{rmk}

Combining this with the upper bound from Prop. \ref{LocalpotentialupperboundII},

\begin{cor}\label{StabilityFermat}
In the Fermat case, there is a uniform $C^0$-stability estimate: 
	\begin{itemize}
		\item
		Inside $U^{s,*}_w\subset X_s$, for $\langle m,w\rangle=1$, the local potentials satisfy
		$
		|\varphi_{CY,s,m} - \psi_{CY,s,m} |\leq  Cs^{-1/2}\log s,
		$
		or equivalently $|\varphi_{CY,s}-\psi_{CY,s}| \leq Cs^{-1/2}\log s$.

		\item
		Inside $U^{s,face}_w$, the local potentials satisfy $|\varphi_{CY,s,0}- \psi_{CY,s,0}| \leq Cs^{-1/2}\log s$, or equivalently $|\varphi_{CY,s}-\psi_{CY,s}|\leq Cs^{-1/2}\log s$.

	\end{itemize}
	
\end{cor}

The point is that in the generic region of $X_s$  \emph{ the Calabi-Yau local potentials are $C^0$-approximated by their regularisations}, which build in convexity by construction, and therefore have a priori Lipschitz bounds.

\section{Fermat case: Metric convergence and SYZ fibration}\label{FermatcasemetricconvergenceandSYZfibration}

We focus on the Fermat family case. We will produce a solution of the real MA equation on $\partial \Delta_\lambda^\vee$ by a subsequential limit, which induces a real MA metric on the regular locus (\cf section \ref{LimitingrealMAmetric}).  Then we show the Calabi-Yau metrics on the degenerating hypersurfaces converge to the real MA metric, both in a $C^\infty_{loc}$-sense (\cf section \ref{Higherregularitysection}) and in the global Gromov-Hausdorff sense (\cf section \ref{GromovHausdorffconvergencesection}). The strong regularity estimates will in particular imply that in the generic region of $X_s$ the CY metrics are collapsing with bounded curvature, which by a result of Zhang \cite{Zhang} allows one to produce a special Lagrangian fibration in the generic region of $X_s$ (\cf section \ref{SpecialLagrangianfibraiongenericregion}).

\subsection{Limiting real MA metric}\label{LimitingrealMAmetric}

We work in the context of section \ref{Localoscillationestimate}, and use the notations therein.
We shall extract some subsequential limit of local potentials for the CY metric $\omega_{CY,s}$, and check that up to a constant it solves the real MA equation on  $\partial \Delta_\lambda^\vee \setminus Sing$  according to Def. \ref{realMAFermat} (\cf also section \ref{RegularitytheoryforrealMA}).

Since the convex functions $u_{CY,s}$ on $N_\R$ produced by double Legendre transform have uniform Lipschitz bounds 
(\ref{Lipschitzboundonu}), by the Arzela-Ascoli theorem we can take a \emph{subsequential limit} as $s\to \infty$, such that $u_{CY,s}\to u_\infty$ in $C^0_{loc}$-topology.  Later we will sometimes suppress mentioning the subsequence for brevity. In particular $u_\infty$ is convex and admissible. We can also pass Cor. \ref{canonicalextensiontranslation} to the limit, to see that in the region $\text{Star}(w)+ \R_{\geq 0} w\subset N_\R$, for any $m$ with $\langle m, w\rangle=1$, the function $u_{\infty,m}=u_\infty-m$ is constant upon translation in the $w$-direction. In particular in such regions the $C^0_{loc}$ convergence improves to $\norm{u_{CY,s,m}- u_{\infty, m} }_{C^0} \to 0$.

By construction 
$
\psi_{CY,s,m}= u_{CY,s,m} \circ \text{Log}_s$, and $u_{CY, s,0}= u_{CY,s}\circ \text{Log}_s.
$
Thus the stability estimate Cor. \ref{StabilityFermat} implies that
\begin{itemize}
\item Inside $U^{s,*}_w\subset X_s$, for $\langle m,w\rangle=1$, the local potentials satisfy
\[
|\varphi_{CY,s,m} - u_{\infty,m}\circ \text{Log}_s |\to 0.
\]

\item
Inside $U^{s,face}_w$, the local potentials satisfy
\[
|\varphi_{CY,s,0}- u_{\infty}\circ \text{Log}_s | \to 0.
\]

\end{itemize}

The rest of this section is devoted to proving

\begin{thm}\label{convergencetorealMA}
On  $\partial\Delta_\lambda^\vee\setminus Sing$,
the locally convex function $u_\infty$ solves the real MA equation in the sense of Def. \ref{realMAFermat} up to a scaling constant: 
\begin{equation}\label{convergencetorealMAequation}
MA(u_\infty)= \frac{ a_\infty }{ \pi^n  n!  } d\mu_\infty,
\end{equation}
where $d\mu_\infty$ is the Lebesgue measure on $\partial\Delta_\lambda^\vee$, and the constant $a_\infty$ is defined by (\ref{ainfty}).
\end{thm}

The intuitive idea is to pass the complex MA equation to some weak limit. The main problem is that the sequence $\varphi_{CY,s}$ live on different manifolds, so we need more effective estimates to pass to the limit.

\begin{lem}\label{pshMArealMA}
Let $u$ be a bounded convex function on the square $\{ |x_i|<1  \}\subset \R^n$. Via the rescaled log map $ s^{-1}\text{Log}: (\C^*)^n\to \R^n$, the function $u$ pulls back to a psh function on $\{  |\log |z_i||<s  \}$. Then the real MA measure of $u$ is related to the pushforward of the complex MA measure of $u\circ s^{-1}\text{Log}$ by
\[
MA(u)= \frac{s^n}{\pi ^n n!}(s^{-1}\text{Log})_*(    \sqrt{-1}\partial\bar{\partial} (u\circ s^{-1}\text{Log})  )^n  . 
\]
\end{lem}

\begin{proof}
If $u$ is smooth, then 
\[
MA(u)(K)= \int_K \det (D^2 u) dx_1\ldots dx_n= \frac{s^n}{\pi ^n n!} \int_{ (s^{-1}\text{Log})^{-1}(K) } ( \sqrt{-1}\partial\bar{\partial} (u\circ s^{-1}\text{Log})  )^n.
\]
Since $u\in C^0$, and both the real and complex MA operators are weakly continuous with respect to $C^0$-limits, this equality passes to general $u$.
\end{proof}

\begin{lem}
(Chern-Levine type estimate)  Let $u$ be a psh function on the annulus region $U=\{ |\log |z_i||<s ,\forall i \} \subset (\C^*)^n$, with $\norm{u}_{L^\infty} \lesssim 1$. Then
\begin{itemize}
\item On the shrinked set $E=\{  |\log |z_i||<s/2    \}$ the measure
\[
\int_E (\sqrt{-1}\partial \bar{\partial} u)^n \leq Cs^{-n}.
\]

\item
Let $u+v$ is another psh function, with $\norm{v}_{L^\infty} \ll 1$. Let $f$ be any compactly supported function on the square $\{ |x_i|<1 \}\subset \R^n$. Then
\[
\int f\{ (\sqrt{-1}\partial \bar{\partial} (u+v)    )^n- (\sqrt{-1}\partial \bar{\partial} u)^n \} \leq Cs^{-n} \norm{f}_{C^2} \norm{v}_{L^\infty}. 
\]
\end{itemize}

\end{lem}

\begin{proof}
Let $\chi$ be a compactly supported nonnegative smooth function on the square $\{ |x_i|<1  \}\subset \R^n$, equal to one on $\{ |x_i|\leq 1/2  \}$. We identify $\chi$ with $\chi\circ s^{-1}\text{Log}$, and denote $\omega_{std}= \sqrt{-1}\sum d\log z_i \wedge d\overline{\log z_i}$. Then
\[
 -Cs^{-2}\omega_{std} \leq \sqrt{-1}\partial \bar{\partial} \chi \leq Cs^{-2}\omega_{std}.
\]
The basic obervation is that if $T$ is a positive current of bidegree $(n-1, n-1)$, then by integration by part,
\[
\begin{split}
&\int_E  \sqrt{-1}\partial \bar{\partial} u \wedge T \leq \int_{\text{supp}(\chi)} \chi \sqrt{-1}\partial \bar{\partial} u \wedge T \\
&=\int_{\text{supp}(\chi)} u \sqrt{-1}\partial \bar{\partial} \chi \wedge T \leq Cs^{-2}\int_{\text{supp}(\chi)}
\omega_{std} \wedge T .
\end{split}
\]
Iterating this argument to lower the power of $\sqrt{-1}\partial\bar{\partial}u$, 
\[
\int_E (\sqrt{-1}\partial \bar{\partial} u)^n \leq Cs^{-2n} \int_U \omega_{std}^n \leq Cs^{-n}.
\]
The second statement is proved similarly by removing $\sqrt{-1}\partial\bar{\partial} v$ factors iteratively.
\end{proof}

\begin{proof}
(Thm. \ref{convergencetorealMA}) There are two subcases: the interior of the top dimensional faces of $\partial \Delta_\lambda^\vee$, and the star of the vertices $\text{Star}(w)$. Since the arguments are almost the same we focus on the latter.

On the interior of $\text{Star}(w)$, we have local affine coordinates $x^{m_1}, \ldots x^{m_n}$, related to the holomorphic $\C^*$-coordinates $z^{m_1},\ldots z^{m_n}$ by $x^{m_i}=s^{-1} \log |z^{m_i}|$. The
star type region $U^{s,*}_w \subset X_s$ can be viewed as a subset of $(\C^*)^n$, so we use the rescaled map $s^{-1} \text{Log}: (\C^*)^n_{z^{m_i}}\to \R^n_{x^{m_i}}$ to pullback the function $u_{\infty,m}$ on $\text{Star}(w)$. On the other hand, $X_s\cap (\C^*)^{n+1}$ maps into $N_\R$ via $\text{Log}_s$, so we can also pullback $u_{\infty,m}$ via $\text{Log}_s$. These two pullbacks differ by at most $Cs^{-1}$ using Cor. \ref{canonicalextensiontranslation}. We also write $\phi=\varphi_{CY,s,m}$.

Take a local test function $f\in C^2_c$ supported in the interior of $\text{Star}(w)$, then $f$ is identified as a local function on $U^{s,*}_w\subset X_s$ via $s^{-1}\text{Log}$. By the Chern-Levine type estimate above,
\[
\begin{split}
&s^n\int f \{  (\sqrt{-1}\partial\bar{\partial} \phi  )^n -  (\sqrt{-1}\partial\bar{\partial} (u_{\infty,m}\circ s^{-1}\text{Log})  )^n  \}
\\
&\leq C \norm{ \phi- u_{\infty,m}\circ s^{-1}\text{Log} }_{L^\infty} \norm{f}_{C^2}\to 0,
\end{split}
\]
as $s\to +\infty$. By the Calabi-Yau condition (\ref{CalabiYaucondition}) and Prop. \ref{canonicalmeasureconvergence}, 
\[
s^n(\sqrt{-1}\partial \bar{\partial} \phi)^n= a_s  d\mu_s = a_\infty   d\mu_s (1+o(1)), \quad s\to +\infty.
\]
Pushing forward via $s^{-1}\text{Log}$, and applying Lemma \ref{pshMArealMA},
\[
\pi^n n!\int f MA(u_\infty)=\lim_{s\to \infty} \int f a_s (s^{-1}\text{Log})_* d\mu_s = a_\infty\int f  d\mu_\infty.
\]
Since this holds for every $f\in C^2_c$, on the interior of this top dimensional face we obtain the measure equality 
 (\ref{convergencetorealMAequation}).
\end{proof}

\subsection{Higher regularity in the generic region}\label{Higherregularitysection}

Once we know the subsequential limit $u_\infty$ satisfies the real MA equation, then by the local regularity theory surveyed in section \ref{RegularitytheoryforrealMA},

\begin{cor}(Regularity of real MA solution)
Inside $\partial \Delta_\lambda^\vee\setminus Sing$, let $\mathcal{R}$ be the set of strictly convex points of $u_\infty$, then $u_\infty\in C^\infty_{loc}(\mathcal{R})$, and the complement of $\mathcal{R}$ is a  closed subset of Hausdorff $(n-1)$-measure zero. In particular $\mathcal{R}$ is path connected, and is open and dense in $\partial \Delta_\lambda^\vee\setminus Sing$.
\end{cor}

\begin{rmk}
In dimension 2, the local regularity theory implies that $\mathcal{R}= \partial \Delta_\lambda^\vee \setminus Sing$, namely the real MA solution is smooth wherever the affine structure is defined. The same might hold in any higher dimension, although this cannot be concluded by local regularity results alone (\cf Remark \ref{Mooneyremark}).
\end{rmk}

We now proceed to a very explicit coordinate version of higher order estimates for the local CY potentials, by transferring regularity from the real MA equation to the complex MA equation.

Let $x\in \mathcal{R}$, then $u_\infty$ (resp. the appropriate $u_{\infty,m})$ has $C^{k,\gamma}$-bound on some coordinate ball $B(x, 2r(x))\subset \mathcal{R}$ contained in a shrinked face (resp. $\text{Star}(w)$). For clarity we focus on the face case. The radius $r(x)$ and the $C^{k,\gamma}$-bound depend on the choice of $x$, but are uniform for $x$ in any fixed compact subset of $\mathcal{R}$.  We identify $u_\infty$ with its pullback to   $(s^{-1}\text{Log})^{-1}(B(x, r(x)))\subset U^{s,face}_w\subset   X_s$.

The local CY potential $ \varphi_{CY,s,0}$  on $(s^{-1}\text{Log})^{-1}(B(x, 2r(x)))$ satisfies
\[
\norm{\varphi_{CY,s,0}- u_\infty\circ s^{-1}\text{Log}  }_{C^0}\to 0, \quad s\to \infty
\]
along the subsequence. We may regard 
 $(s^{-1}\text{Log})^{-1}(B(x, 2r(x)))$ as an open subset of $(\C^*)^n$. On the universal cover of $(\C^*)^n$, we use the natural coordinates $s^{-1} \log z^{m_i}$ for $i=1,\ldots n$.

Now $\varphi_{CY,s,0}$ satisfies the complex MA equation (\cf (\ref{CalabiYaucondition})(\ref{normalisedcanonicalmeasure}))
\[
(\sqrt{-1}\partial\bar{\partial} \varphi_{CY,s,0})^n= a_s s^{-n} d\mu_s= \frac{ a_s }{(4\pi s^2)^n  } \sqrt{-1}^{n^2} \Omega_s\wedge \overline{\Omega_s}.
\]
By the holomorphic volume form formula (\ref{holomorphicvolumetoricregion}),
\[
(\sqrt{-1}\partial\bar{\partial}  \varphi_{CY,s,0})^n= \frac{ a_s }{(4\pi )^n  }  (1+o(1)) \prod_i \sqrt{-1} s^{-1}d\log z^{m_i} \wedge s^{-1}d\overline{\log z^{m_i}},
\]
where the $o(1)$ term in fact  has exponentially small $C^\infty$ bounds in $s^{-1}\log z^{m_i}$ coordinates; the higher order bound uses that $\Omega_s$ is holomorphic. On the other hand by the calculation in section \ref{LimitingrealMAmetric}, the pullback of $u_\infty$ satisfies
\[
(\sqrt{-1}\partial\bar{\partial} (u_\infty \circ s^{-1}\text{Log} ) )^n= \frac{ a_\infty }{ (4\pi)^n  }  \prod_i s^{-1}\sqrt{-1}d\log z^{m_i} \wedge s^{-1}d\overline{\log z^{m_i}}.
\]
To summarize, the deviation of RHS is negligible and the deviation between $\varphi_{CY,s,0}$ and $u_\infty\circ s^{-1}\text{Log}$ is small in $C^0$-norm. Applying Savin's Thm. \ref{Savin}, 

\begin{thm}\label{smoothconvergence}(Smooth convergence in generic regions)
As $s\to +\infty$ along the subsequence, assume the coordinate ball $B(x, 2r(x))\subset \mathcal{R}$, then on the region $(s^{-1}\text{Log})^{-1}(B(x,r(x)))\subset X_s$, we have the following higher regularity estimates with respect to the $C^{k,\gamma}$-norm in the $s^{-1}\log z^{m_i}$ coordinates.
\begin{itemize}
\item  In the face type region $U_w^{s,face}$ case \[
\norm{ \varphi_{CY,s,0}-u_\infty\circ s^{-1}\text{Log} }_{C^{k,\gamma}( (s^{-1}\text{Log})^{-1}(B(x,r(x) ) ) } \to 0.
\]

\item  In the star type region $U_w^{s,*}$ case, for $\langle m, w\rangle =1$,
\[
\norm{ \varphi_{CY,s,m}-u_{\infty,m}\circ s^{-1}\text{Log} }_{C^{k,\gamma}( (s^{-1}\text{Log})^{-1}(B(x,r(x) )   )} \to 0.
\]
\end{itemize}
 The convergence rate is uniform for $x$ on any fixed compact subset of $\mathcal{R}$. 
\end{thm}

The intuition is that in the generic regular locus in the toric part of  $X_s$, the local CY potentials converge in some $C^\infty_{loc}$ sense. 

\begin{Notation}
	For every compact $K\subset \mathcal{R}$, let $U_{s,K}$ denote the union of the regions $(s^{-1}\text{Log})^{-1} (B(x,r(x))  )$ for $x\in K$; the convergence rates will be uniform on $U_{s,K}$. Notice that 
	\[
	\limsup_{s\to \infty }
	\frac{ \text{Vol}( U_{s,K} ) }  { \text{Vol}(X_s)  } 
	\geq  \frac{ \int_K d\mu_\infty } { \int_{ \partial \Delta_\lambda^\vee} d\mu_\infty } , 
	\]
	so by taking a compact exhaustion of $\mathcal{R}$, we may assume  $U_{s,K}$ occupies a percentage of the total measure arbitrarily close to 1.
\end{Notation}

\begin{rmk}
If one can show that the limiting real MA metric is unique, then there will be no need to pass to a subsequence.
\end{rmk}

Next we discuss CY metrics in $(s^{-1}\text{Log})^{-1}(B(x,r(x)))\subset U_{s,K}$.
\begin{itemize}
\item In the face type region case, up to $C^\infty$-small error in the $s^{-1}\log z^{m_i}$ coordinates,
\[
\begin{split}
& \omega_{CY,s}=\sqrt{-1}\partial \bar{\partial} \varphi_{CY,s,0}
\\
\approx &\sqrt{-1}\partial \bar{\partial} u_\infty \circ s^{-1}\text{Log}= \frac{1}{4} \frac{ \partial^2 u_\infty }{ \partial x^{m_i} \partial x^{m_j}   }\sqrt{-1} s^{-1}d\log z^{m_i} \wedge s^{-1} d\overline{ \log z^{m_j} },
\end{split}
\]
hence the CY metrics $g_{CY,s}$ is up to $C^\infty$-small error
\begin{equation}\label{CYmetricsemiflat1}
g_{CY,s}\approx \text{Re} \{ \frac{1}{2} \frac{ \partial^2 u_\infty }{ \partial x^{m_i} \partial x^{m_j}   }   s^{-1}d\log z^{m_i} \otimes s^{-1} d\overline{ \log z^{m_j} }  \}.
\end{equation}

\item
Likewise in the star type region case, up to $C^\infty$ small error in the $s^{-1}\log z^{m_i}$ coordinates,
\begin{equation}\label{CYmetricsemiflat2}
\begin{cases}
\omega_{CY,s} \approx \frac{1}{4} \frac{ \partial^2 u_{\infty,m} }{ \partial x^{m_i} \partial x^{m_j}   }\sqrt{-1} s^{-1}d\log z^{m_i} \wedge s^{-1} d\overline{ \log z^{m_j} },
\\
g_{CY,s}\approx \text{Re} \{ \frac{1}{2} \frac{ \partial^2 u_{\infty,m} }{ \partial x^{m_i} \partial x^{m_j}   }   s^{-1}d\log z^{m_i} \otimes s^{-1} d\overline{ \log z^{m_j} }  \}.
\end{cases}
\end{equation}
\end{itemize}

Notice in such local $(\C^*)^n$ coordinates, the rescaled log map $s^{-1}\text{Log}$  gives a local $T^n$-fibration. The metric associated to $\sqrt{-1}\partial \bar{\partial} (u_\infty\circ s^{-1}\text{Log}) $ is a \emph{semiflat metric}, namely a $T^n$-invariant metric which is flat when restricted to any $T^n$-fibre. Thus (\ref{CYmetricsemiflat1})(\ref{CYmetricsemiflat2}) assert that \emph{the Calabi-Yau metrics $g_{CY,s}$ are $C^\infty$-approximated by semiflat metrics in the regular regions}.

\begin{cor}\label{boundedsectionalcurvature}
On $U_{s,K}\subset X_s$ the sectional curvature has a uniform bound $|\text{Riem}(g_{CY,s})|\leq C$, and the injectivity radius  satisfies $C^{-1}s^{-1} \leq \text{inj} \leq Cs^{-1}$, with constants depending on $K\subset \mathcal{R}$.
\end{cor}



\subsection{Gromov-Hausdorff convergence}\label{GromovHausdorffconvergencesection}

On the regular locus $\mathcal{R}\subset \partial \Delta_\lambda^\vee$ we have a well defined real MA metric, 
\begin{equation}
g_\infty= 
\begin{cases}
\frac{1}{2}\sum_{i,j}\frac{\partial^2 u_\infty }{\partial x_i \partial x_j} dx_i dx_j, \quad \text{ on the face regions},
\\
\frac{1}{2}\sum_{i,j}\frac{\partial^2 u_{\infty,m} }{\partial x_i \partial x_j} dx_i dx_j, \quad \text{ on the star regions}.
\end{cases}
\end{equation}
Notice the definitions are compatible on overlapping regions.
Let $(\bar{\mathcal{R}}, g_\infty)$ be the metric completion. The metric asymptotes  (\ref{CYmetricsemiflat1})(\ref{CYmetricsemiflat2}) say that in some $C^\infty_{loc}$ sense the collapsing CY metrics $g_{CY,s}$ converge to the metric $g_\infty$ on $\mathcal{R}$, and we know $\mathcal{R}$ is path connected because its complement has zero $\mathcal{H}^{n-1}$-measure.

\begin{rmk}
We do not know if $\bar{\mathcal{R}}$ is homeomorphic to $\partial \Delta_\lambda^\vee\simeq S^n$, as the regularity theory of the real MA equation on a singular affine manifold is not yet developed, and we know little about what can happen near singularities.
\end{rmk}

The goal of this section is to show

\begin{thm}\label{GHconvergencethm}
The subsequence of collapsing CY metrics $(X_s, g_{CY,s})$ converges in the Gromov-Hausdorff sense to $(\bar{\mathcal{R}}, g_\infty)$.
\end{thm}

\begin{prop}\label{diameterbound}
There is a uniform diameter bound 
\[
\text{diam}(X_s, g_{CY,s}) \leq C.
\]
\end{prop}

\begin{proof}
This argument is essentially the same as \cite[Thm 3.1]{Tosattidiameter}. 
We quote \cite[Lem 3.2]{Tosattidiameter}:
\begin{lem}
Let $(M^{2n}, g)$ be a closed Riemannian manifold with $ Ric(g) \geq 0$, let $ p \in M $and $1 < R \leq diam(X, g)$. Then
$ \frac{R -1}{4n}
\leq
\frac{ \text{Vol}(B(p, 2(R + 1))) }
{ \text{Vol}(B(p, 1))}$.	
\end{lem}

Using Thm. \ref{smoothconvergence}, we can find inside the regular region of $X_s$ some geodesic ball $B_{g_{CY,s}}(p,r)$ of radius $r<1$, occupying a nontrivial portion of the total volume:
\[
\frac{ \text{Vol}(B_{g_{CY,s}}(p, r))) }
{ \text{Vol}(X_s)} \geq \epsilon>0,
\]
with $\epsilon$ independent of $s$. Now applying the Lemma to the rescaled  CY metric $r^{-2}g_{CY,s}$, 
\[
\frac{ \text{diam}(X_s) -r}{4n r}
\leq
\frac{ \text{Vol}(B_{g_{CY,s}}(p, 2(\text{diam}(X_s) + r))) }
{ \text{Vol}(B_{g_{CY,s}}(p, r))} \leq \frac{ \text{Vol}(X_s) }
{ \text{Vol}(B_{g_{CY,s}}(p, r))} \leq \epsilon^{-1},
\]
so $\text{diam}(X_s) \leq Cr \leq C$ as required.
\end{proof}

\begin{proof}
(Thm. \ref{GHconvergencethm}) By Thm \ref{smoothconvergence} we already know the metric convergence over any properly contained open subset of $\mathcal{R}$, which corresponds to a region $U_s\subset X_s$, with nearly the full measure:
\[
\text{Vol}(U_s)>(1-\epsilon) \text{Vol}(X_s),
\]
where $\epsilon$ can be chosen arbitrarily small. It now suffices to show any point $p\in X_s\setminus U_s$ is close to $U_s$. For any $r>0$ such that the geodesic ball $B_{g_{CY,s}}(p,r) \subset X_s\setminus U_s  $, the Bishop-Gromov inequality implies
\[
\left(\frac{r}{ \text{diam}(X_s) } \right)^{2n} \leq
\frac{ \text{Vol}(B_{g_{CY,s} }(p,r) ) }{  \text{Vol}(X_s)  } \leq \frac{ \text{Vol}(X_s\setminus U_s ) }{  \text{Vol}(X_s)  }<\epsilon.
\]
Taking the sup of all such $r$,
\[
\text{dist}_{g_{CY,s}}(p, U_s) \leq \epsilon^{1/2n} \text{diam}(X_s) \leq C \epsilon^{1/2n},
\]
which can be made arbitrarily small.
\end{proof}

\subsection{Special Lagrangian fibration in the generic region}\label{SpecialLagrangianfibraiongenericregion}

In the setting of section \ref{Higherregularitysection}, the very strong regularity bounds in the generic region leads to the existence of special Lagrangian $T^n$-fibrations thereon.

\begin{thm}
For any fixed compact $K\subset \mathcal{R}$, for $s\gg 1$ depending on $K$, there is a special Lagrangian (SLag) $T^n$-fibration on an open subset of $X_s$ containing $U_{s,K}$.
\end{thm}

\begin{rmk}
By considering a compact exhaustion of $\mathcal{R}$, we can choose $K$ so that the region $U_{s,K}$ occupies a percentage of the total measure on $X_s$ arbitrarily close to 1.
\end{rmk}

\begin{proof}
Since $K$ is a compact subset in the open set $\mathcal{R}$, we can find an open set $\mathcal{U}\subset K$ properly contained in $\mathcal{R}$. This ensures that the smooth convergence in Thm. \ref{smoothconvergence} happens uniformly on a slightly larger set $U_{s, K' }$ than $U_{s,K}$. We  assume $s\gg 1$ as ususal.

Consider a coordinate region $(s^{-1}\text{Log})^{-1}(B(x,r(x))$ contained in this larger set, which is topologically $T^n\times B(x, r(x))$. Here the $T^n$ is well defined as a homology cycle independent of the coordinates. We define the phase angles $\theta_s$ by requiring
$
\int_{T^n} e^{\sqrt{-1} \theta_s } \Omega>0.
$
We consider the rescaled CY metrics $(s^2 g_{CY,s}, s^2 \omega_{CY,s})$, so the diameter of $T^n$ fibres are now of order $O(1)$ by (\ref{CYmetricsemiflat1})(\ref{CYmetricsemiflat2}). Within any log scale, these rescaled CY structures are $C^\infty$-close to the standard flat structures in section \ref{SpecialLagrangiansurvey} up to constant factors. By construction the K\"ahler forms are exact in these coordinate charts. Thus by Zhang's result surveyed in section \ref{SpecialLagrangiansurvey}, within any log scale, we can construct a SLag $T^n$-fibration with phase $\theta_s$, whose fibres are very small $C^\infty$-perturbations of the fibres of the map $(\C^*)^n\to \R^n$,
\[
\text{Log}: (z^{m_1}, \ldots , z^{m_n}) \to ( \log |z^{m_1}|, \ldots \log |z^{m_n}|   ).
\]

Observe that on overlapping charts, the $\text{Log}$-fibres with respect to one chart are very small $C^\infty$-perturbations of the $\text{Log}$-fibres of the other chart. Then the uniqueness part of Zhang's argument shows that on overlapping charts the SLag $T^n$-fibrations are in fact defined independent of charts. (It is the local universal family of SLags within the perturbative regime.) Thus the local constructions glue to a SLag fibration on a subset of $X_s$ containing $U_{s,K}$ as required.
\end{proof}


\begin{thebibliography}{7}
	



\bibitem{Blocki} Błocki, Zbigniew; Kołodziej, Sławomir. On regularization of plurisubharmonic functions on manifolds. Proc. Amer. Math. Soc. 135 (2007), no. 7, 2089--2093.






\bibitem{Blockilecture} 
Błocki, Zbigniew. The Calabi-Yau theorem. Complex Monge-Ampère equations and geodesics in the space of Kähler metrics, 201--227, Lecture Notes in Math., 2038, Springer, Heidelberg, 2012.




\bibitem{Boucksom1} 
Boucksom, Sébastien; Jonsson, Mattias. Tropical and non-Archimedean limits of degenerating families of volume forms. J. Éc. polytech. Math. 4 (2017), 87--139. 




\bibitem{Boucksom} 
Boucksom, Sébastien; Favre, Charles; Jonsson, Mattias. Solution to a non-Archimedean Monge-Ampère equation. J. Amer. Math. Soc. 28 (2015), no. 3, 617--667.


\bibitem{Boucksomsurvey}
Boucksom, Sébastien; Favre, Charles; Jonsson, Mattias. The non-Archimedean Monge-Ampère equation. Nonarchimedean and tropical geometry, 31--49, Simons Symp., Springer, [Cham], 2016. 


\bibitem{Caff1}  Caffarelli, L. A. A localization property of viscosity solutions to the Monge-Ampère equation and their strict convexity. Ann. of Math. (2) 131 (1990), no. 1, 129--134. 

\bibitem{Caff2} Caffarelli, Luis A. Interior $W^{2,p}$ estimates for solutions of the Monge-Ampère equation. Ann. of Math. (2) 131 (1990), no. 1, 135--150. 





\bibitem{Caff4} Caffarelli, Luis A. A note on the degeneracy of convex solutions to Monge Ampère equation. Comm. Partial Differential Equations 18 (1993), no. 7-8, 1213--1217. 



\bibitem{CaffarelliViaclovsky} 
Caffarelli, Luis A.; Viaclovsky, Jeff A. On the regularity of solutions to Monge-Ampère equations on Hessian manifolds. Comm. Partial Differential Equations 26 (2001), no. 11-12, 2339--2351.



\bibitem{CollinsTosatti} Collins, Tristan C.; Tosatti, Valentino. An extension theorem for Kähler currents with analytic singularities. Ann. Fac. Sci. Toulouse Math. (6) 23 (2014), no. 4, 893--905.






\bibitem{Coman} Coman, Dan; Guedj, Vincent; Zeriahi, Ahmed. Extension of plurisubharmonic functions with growth control. J. Reine Angew. Math. 676 (2013), 33--49. 




\bibitem{DemaillyPali}
Demailly, Jean-Pierre; Pali, Nefton. Degenerate complex Monge-Ampère equations over compact Kähler manifolds. Internat. J. Math. 21 (2010), no. 3, 357--405. 




\bibitem{Donaldsontoric} Donaldson, Simon K. Kähler geometry on toric manifolds, and some other manifolds with large symmetry. Handbook of geometric analysis. No. 1, 29--75, Adv. Lect. Math. (ALM), 7, Int. Press, Somerville, MA, 2008.






\bibitem{DonaldsonSun} 

Donaldson, Simon; Sun, Song. Gromov-Hausdorff limits of Kähler manifolds and algebraic geometry. Acta Math. 213 (2014), no. 1, 63--106.



\bibitem{EGZ} 
Eyssidieux, Philippe; Guedj, Vincent; Zeriahi, Ahmed. Singular Kähler-Einstein metrics. J. Amer. Math. Soc. 22 (2009), no. 3, 607--639.



\bibitem{EGZ2}
Eyssidieux, Philippe; Guedj, Vincent; Zeriahi, Ahmed. A priori $L^\infty$-estimates for degenerate complex Monge-Ampère equations. Int. Math. Res. Not. IMRN 2008, Art. ID rnn 070, 8 pp.


\bibitem{Lorenzo} 
Foscolo, Lorenzo. ALF gravitational instantons and collapsing Ricci-flat metrics on the $K3$ surface. J. Differential Geom. 112 (2019), no. 1, 79--120.



\bibitem{Gross} 
Gross, Mark. Mirror symmetry and the Strominger-Yau-Zaslow conjecture. Current developments in mathematics 2012, 133--191, Int. Press, Somerville, MA, 2013. 



\bibitem{Gross2} 
Gross, Mark. Topological mirror symmetry. Invent. Math. 144 (2001), no. 1, 75--137. 




\bibitem{GrossTosattiZhang} 
Gross, Mark; Tosatti, Valentino; Zhang, Yuguang. Collapsing of abelian fibered Calabi-Yau manifolds. Duke Math. J. 162 (2013), no. 3, 517--551. 




\bibitem{GrossWilson}
Gross, Mark; Wilson, P. M. H. Large complex structure limits of $K3$ surfaces. J. Differential Geom. 55 (2000), no. 3, 475--546.



\bibitem{GZ} 
Guedj, Vincent; Zeriahi, Ahmed. Intrinsic capacities on compact Kähler manifolds. J. Geom. Anal. 15 (2005), no. 4, 607--639.


\bibitem{HaaseZharkov} 
Haase, Christian; Zharkov, Ilia. Integral affine structures on spheres and torus fibrations of Calabi-Yau toric hypersurfaces I. 	arXiv:math/0205321.

 

\bibitem{HaaseZharkov2} Haase Zharkov 2
Haase, Christian; Zharkov, Ilia.
Integral affine structures on spheres and torus fibrations of Calabi-Yau toric hypersurfaces II. 	arXiv:math/0301222.



\bibitem{HarveyLawson} 
Harvey, Reese; Lawson, H. Blaine, Jr. Calibrated geometries. Acta Math. 148 (1982), 47--157. 




\bibitem{HSVZ}
Hein, Hans-Joachim;    Sun, Song;  Viaclovsky, Jeff;  Zhang, Ruobing. Nilpotent structures and collapsing Ricci-flat metrics on K3 surfaces. 	arXiv:1807.09367. 







\bibitem{Joyce} 
Joyce, Dominic. Singularities of special Lagrangian fibrations and the SYZ conjecture. Comm. Anal. Geom. 11 (2003), no. 5, 859--907.

\bibitem{KS2}
Kontsevich, Maxim; Soibelman, Yan. Homological mirror symmetry and torus fibrations. Symplectic geometry and mirror symmetry (Seoul, 2000), 203--263, World Sci. Publ., River Edge, NJ, 2001.



\bibitem{KS} 
Kontsevich, Maxim; Soibelman, Yan. Affine structures and non-Archimedean analytic spaces. The unity of mathematics, 321--385, Progr. Math., 244, Birkhäuser Boston, Boston, MA, 2006. 



\bibitem{Li} Li, Yang.
SYZ geometry for Calabi-Yau 3-folds: Taub-NUT and Ooguri-Vafa type metrics. 
arXiv:1902.08770.


\bibitem{Mooney} 
Mooney, Connor. Partial regularity for singular solutions to the Monge-Ampère equation. Comm. Pure Appl. Math. 68 (2015), no. 6, 1066--1084.





\bibitem{Odaka} 
Odaka, Yuji; Oshima, Yoshiki. Collapsing K3 surfaces and Moduli compactification. Proc. Japan Acad. Ser. A Math. Sci. 94 (2018), no. 8, 81--86.


\bibitem{Savin} 
Savin, Ovidiu. Small perturbation solutions for elliptic equations. Comm. Partial Differential Equations 32 (2007), no. 4-6, 557--578.



\bibitem{SYZ} 
Strominger, Andrew; Yau, Shing-Tung; Zaslow, Eric. Mirror symmetry is $T$-duality.
Nucl.Phys.B479:243-259,1996.




\bibitem{HSVZ2} 
Sun, Song; Zhang, Ruobing. Complex structure degenerations and collapsing of Calabi-Yau metrics. 	arXiv:1906.03368.

\bibitem{Tosattidiameter} 
Tosatti, Valentino. Limits of Calabi-Yau metrics when the Kähler class degenerates. J. Eur. Math. Soc. (JEMS) 11 (2009), no. 4, 755--776.

\bibitem{Tosatti}
Tosatti, Valentino. Adiabatic limits of Ricci-flat Kähler metrics. J. Differential Geom. 84 (2010), no. 2, 427--453.








\bibitem{Yau}
Yau, Shing Tung. On the Ricci curvature of a compact Kähler manifold and the complex Monge-Ampère equation. I. Comm. Pure Appl. Math. 31 (1978), no. 3, 339--411.




\bibitem{Zharkov} 
Zharkov, Ilia. Limiting behavior of local Calabi-Yau metrics. Adv. Theor. Math. Phys. 8 (2004), no. 3, 395--420. 

\bibitem{Zeriahi}
Zeriahi, Ahmed. Volume and capacity of sublevel sets of a Lelong class of plurisubharmonic functions. Indiana Univ. Math. J. 50 (2001), no. 1, 671--703. 



\bibitem{Zhang}
Zhang, Yuguang. Collapsing of Calabi-Yau manifolds and special Lagrangian submanifolds. Univ. Iagel. Acta Math. No. 54 (2017), 53--78.




\end{thebibliography}
\end{document}